\documentclass [a4paper,twoside,12pt]{article}

\usepackage[utf8]{inputenc}
\usepackage[T1]{fontenc}
\usepackage{amsmath,amsfonts,amssymb}
\usepackage{vmargin,graphicx,theorem}
\usepackage[english]{babel}
\usepackage{enumerate}
\usepackage{color}
\usepackage{pstricks}
\usepackage{mdwlist}
\RequirePackage[colorlinks,linkcolor=blue,citecolor=red,urlcolor=blue]{hyperref}

\setpapersize[portrait]{A4}
\setmarginsrb{1.5cm}{1cm}{1.5cm}{2.5cm}{1.5cm}{0cm}{0.5cm}{2cm}
\newcommand{\disp}{\displaystyle}
\newcommand{\dL}{\ensuremath{\mathbb{L}}}

\newcommand{\dR}{\ensuremath{\mathbb{R}}}

\newtheorem{ethm}{Theorem}%[section]

\newtheorem{eprop}[ethm]{Proposition}

\newtheorem{elem}[ethm]{Lemma}

\newtheorem{edefi}[ethm]{Definition}

\newtheorem{erem}[ethm]{Remark}

\newtheorem{eex}[ethm]{Example}

\newtheorem{eass}{Regularity assumption}

\newcommand{\proofend}{~$\rhd$}
\newcommand{\proofbegin}{~$\lhd$}

\newenvironment{Proof}{\par\begin{trivlist}%
\item[]{\bf Proof.\ }}%
{\hfill $\square$ \end{trivlist}\par}
\newenvironment{tProof}[1]{\par\begin{trivlist}%
\item[]{\bf Proof of #1.\ }}%
{\hfill $\square$ \end{trivlist}\par}

\newcommand{\p}[4]{{#3}\!\left#1{#4}\right#2}

\newcommand{\PAR}[1]{\ensuremath{{\left(#1\right)}}} % (1)
\newcommand{\SBRA}[1]{\ensuremath{{\left[#1\right]}}} % [1]
\newcommand{\BRA}[1]{\ensuremath{{\left\{#1\right\}}}} % {1}
\renewcommand{\phi}{\varphi}
\newcommand{\ep}{{\varepsilon}} % epsilon

\renewcommand{\geq}{\geqslant}

\newcommand{\entf}[1]{{\rm{Ent}}_{#1}}

\newcommand{\ent}[2]{\p(){\entf{#1}}{#2}}

\def\disp{\displaystyle}

\newcommand{\N}{\ensuremath{\mathbb{N}}}
\newcommand{\M}{\ensuremath{\mathbf{M}}}
\newcommand{\X}{\ensuremath{\mathbf{X}}}
\newcommand{\R}{\dR}

\newcommand{\e}{\varepsilon}

\newcommand{\beq}{\begin{equation}}\newcommand{\eeq}{\end{equation}}
\begin{document}

\title{Equivalence between dimensional contractions in Wasserstein distance and the curvature-dimension condition. }
\author{ Fran\c{c}ois Bolley\thanks{Laboratoire de probabilit\'es et mod\`eles al\'eatoires, Umr Cnrs 7599, Univ. Pierre et Marie Curie - Paris 6, France. francois.bolley@upmc.fr}, Ivan Gentil\thanks{Institut Camille Jordan, Umr Cnrs 5208, Univ. Claude Bernard Lyon 1, France. gentil@math.univ-lyon1.fr}, Arnaud Guillin\thanks{Institut Universitaire de France and Laboratoire de Math\'ematiques, Umr Cnrs 6620, Univ. Blaise Pascal, France. guillin@math.univ-bpclermont.fr}, Kazumasa Kuwada\thanks{Tokyo Institute of Technology, Japan. kuwada@math.titech.ac.jp}}

\date{\today}

\maketitle

\abstract{
The curvature-dimension condition is a generalization of the Bochner inequality to weighted Riemannian manifolds and general metric measure spaces. It is now known 
to be equivalent to evolution variational inequalities for the heat semigroup, and quadratic Wasserstein distance contraction properties at different times.
On the other hand, in a compact Riemannian manifold, it implies a same-time Wasserstein contraction property for this semigroup. 
In this work we generalize the latter result to metric measure spaces and more importantly  prove the converse: contraction inequalities are equivalent to curvature-dimension conditions.
Links with functional inequalities are also investigated.
}
\bigskip

\bigskip

\section*{Introduction}

The von Renesse-Sturm theorem (see \cite{sturm-vonrenesse}) ensures that a Wasserstein distance contraction property between solutions to the heat equation on a Riemannian manifold is equivalent to a lower curvature condition. This result is one of the first equivalence results relating the Wasserstein distance and a curvature condition.  Recent works have been devoted to a more precise curvature-dimension condition instead of a sole curvature condition.  In this work, and  in a fairly general framework, we derive new {\it dimensional} contraction properties under a curvature-dimension condition and we show that they are all equivalent to it.

 \medskip

  Let  $\Delta$ be the Laplace-Beltrami operator on a smooth Riemannian manifold $(\M, \mathcal G)$ and let $(P_t f)_{t \geq 0}$ be the solution to the heat equation $\partial_tu=\Delta u$ with $f$ as the initial condition. Many of the coming notions and results have been considered in a more general setting, but for simplicity in the introduction we focus on this case. The Bochner identity states that
$$
\frac{1}{2}\Delta|\nabla f|^2-\nabla f\cdot\nabla \Delta f=|\nabla\nabla f|^2+\mathrm{Ric}(\nabla f,\nabla f)
$$
where $\mathrm{Ric}$ is the Ricci curvature of $(\M,\mathcal G)$. The manifold associated with its Laplacian is said to satisfy the $CD(R,m)$  curvature-dimension condition if its Ricci curvature is uniformly bounded from below by $R \in \mathbb R$ and its dimension is smaller than $m  \in (0, + \infty]$. In this case
\begin{equation}
\label{eq-bochner}
\frac{1}{2}\Delta|\nabla f|^2-\nabla f\cdot\nabla \Delta f\geq \frac1m(\Delta f)^2+R|\nabla f|^2
\end{equation}
by the Cauchy-Schwarz inequality.
 The $CD(R,m)$ condition and~\eqref{eq-bochner} are the starting point of many comparison theorems, functional and geometrical inequalities, bounds on the heat kernel, etc. (see e.g.~\cite{bgl-book,EKS13,villani-book2,wang-book}). 

\medskip

 In this work we focus on the link between the curvature-dimension condition and Wasserstein distance contraction properties of the heat semigroup.
The von Renesse-Sturm theorem \cite{sturm-vonrenesse} states that: the $CD(R,\infty)$ condition holds if and only if
\begin{equation}
\label{eq-premiere}
 W_2^2(P_tfdx,P_tgdx)\leq e^{-2Rt}W_2^2(fdx,gdx)
\end{equation}
for all $t \geq 0$ and probability densities  $f,g$ with respect to the Riemannian measure $dx$. Here $W_2$ is the Wasserstein distance with quadratic cost.
 
There are many proofs of this result as well as extensions to more general evolutions and spaces, see for instance~\cite{ambrosio-gigli-savare,bgl-15,bgl-book,gko13,kuwada10,otto05,wang-book,wang11}. Following the seminal papers~\cite{LV,sturm}, attention has been drawn to taking the {\it dimension} of the manifold into account. 

A first way of including the dimension is to use {\it two different times} $s$ and $t$ in the inequality~\eqref{eq-premiere}. It is proved in~\cite{bgl-15,kuwada15} that the $CD(0,m)$ condition implies 
\begin{equation}
\label{eq-cas-simple}
 W_2^2(P_sfdx,P_tgdx)\leq W_2^2(fdx,gdx)+2m(\sqrt{t}-\sqrt{s})^2
\end{equation}
for all $s,t\geq0$ and all probability densities $f,g$. 
A non zero lower bound on the curvature and the equivalence have been further considered in~\cite{EKS13,kuwada15}:

\begin{itemize}
\item In~\cite{kuwada15}, the fourth author proved that the $CD(R,m)$ condition holds  if and only if 
\begin{equation}
\label{eq-last}
W_2^2(P_t fdx,P_sgdx)\leq A(s,t,R,m) W_2^2( fdx,gdx)+B(s,t,m,R)
\end{equation}
for all $s,t\geq0$ and all probability densities $f,g$, and for appropriate positive functions $A, B.$ 

\item In~\cite{EKS13}, the authors proved that the $CD(R,m)$ condition holds  if and only if
\begin{multline}
\label{eq-eks}
s_{\frac Rm}\left(\frac 12 W_2(P_tf dx,P_sg dx)\right)^2
\leq e^{-R(t+s)}\,s_{\frac Rm}\left(\frac 12 W_2(fdx,gdx)\right)^2
\\
+\frac mR(1-e^{-R(s+t)})\frac{(\sqrt{t}-\sqrt{s})^2}{2(t+s)}
\end{multline}
for all $s,t\geq0$ and all probability densities $f, g$. Here $s_r(x)=\sin(\sqrt{r}x) / \sqrt{r}$ if $r>0$, $s_r(x)=\sinh(\sqrt{\vert r \vert}x) / \sqrt{\vert r \vert}$ if $r<0$ and $s_0(x)=x$, hence recovering~\eqref{eq-cas-simple} when $R=0$. Both inequalities~\eqref{eq-last} and~\eqref{eq-eks} are extensions of~\eqref{eq-premiere} and~\eqref{eq-cas-simple}, taking the dimension into account.
 \end{itemize}

\medskip
Contraction properties with the {\it same time} have been derived in~\cite{bgg28} for the Euclidean heat equation in $\R^m$,  and then extended by the third author in~\cite{G15} to a compact Riemannian manifold. Let $\ent{dx}{h} = \int h \, \log h \, dx$ be the entropy of a probability density $h$. Then the $CD(R,m)$ condition implies
$$
W_2^2(P_t fdx,P_t gdx)\leq e^{-2R t}\,W_2^2(fdx,gdx)\\-\frac{2}{m}\int_0^t \! e^{-2R(t-u)}\PAR{\ent{dx}{P_u g}-\ent{dx}{P_u f}}^2du
$$
for all $t \geq 0$ and all $f$, $g$ probability densities.
This bound has also been proved in~\cite{bgg28} for the Markov transportation distance instead of the $W_2$ distance. This distance differs from $W_2$ and has actually been  tailored to Markov semigroups and the Bakry-\'Emery $\Gamma_2$ calculus. Dimensional contraction properties for a Wasserstein distance defined with an adapted cost have also been derived in~\cite{wang11}.

\bigskip

In this paper we derive diverse  {\it same time} contraction inequalities under a general $CD(R,m)$ curvature-dimension condition, and in fact prove that they are all {\it equivalent} to this condition. The results and the proof will be given in the two settings of a smooth Riemannian manifold and of a more general metric measure space, more precisely in the setting introduced in~\cite{AGS_BE} of a Riemannian energy measure space. 

The paper is organized as follows. In Section~\ref{sec-main-result}, we state and explain the context of  our main result, Theorem~\ref{thm-legros}. In Section~\ref{sec-proofdebut}, we prove the easier implications, leaving the main issue aside: from the weakest contraction to the curvature-dimension condition. Some arguments in this section require a detailed formulation given in Section~\ref{sec-mms} below : thus they are only outlined there and complemented in Section~\ref{subsec:3-4}. In Section~\ref{sec-strategy}, we present the strategy of our proof, motivated by the elementary gradient flow approach in Euclidean space. The result is proved on a Riemannian manifold in Section~\ref{sec-riemannian}, and on a Riemannian energy measure space in Section~\ref{sec-mms}. The general strategy is the same in both settings, and it could seem redundant to give both proofs. However the proof in the Riemannian setting is rather simpler, presents the most important steps of the argument and thus gives a way to get it in a more general space.  We believe that it is an opportunity to emphasize, in our example, the main issues arising in transferring a proof in the Riemannian setting to the abstract measure space setting. Indeed, there, regularity is no more available ``for free'', and our proof will crucially use a whole panel of powerful tools developed by L.~Ambrosio, N.~Gigli, G.~Savar\'e, K.-T.~Sturm and coauthors to overcome this difficulty, in particular localization and mollification by semigroup.
The last section gives a new and simple derivation of a classical entropy-energy inequality, as well as dimensional HWI inequalities: for this we start from our contraction inequalities instead of the curvature-dimension condition, as in earlier works.

\section{Main result}
\label{sec-main-result}

Our main theorem states that, in a quite general framework, a curvature-dimension condition is equivalent to same time Wasserstein distance contraction inequalities. 

\medskip

Let $(\X,d)$ be a Polish metric space, $\mathcal P(\X)$ be the set of Borel probability measures on $\X$ and 
$\mathcal P_2(\X)$ be the set of all $\mu\in\mathcal P(\X)$ such that $\int d(x_0,x)^2 \, d\mu(x)<\infty$ for some $x_0\in \X.$ 
The (quadratic) Wasserstein distance between $\nu_1$ and $\nu_2$ in $\mathcal P_2(\X)$ is defined by  
$$
W_2(\nu_1,\nu_2)=\inf_{\pi} \sqrt{\iint d(x,y)^2 \, d\pi(x,y)}
$$
where the infimum runs over all probability measures $\pi$ on $\X\times \X$ with marginals $\nu_1$ and $\nu_2$. 

A fundamental tool is the Kantorovich dual representation : 
for $\nu_1,\nu_2\in\mathcal P_2(\X)$, 
\begin{equation}
\label{eq-kanto}
\frac{W_2^2(\nu_1,\nu_2)}{2}=\sup_{\psi} \Big\{ \int Q\psi \, d\nu_1-\int\psi \, d\nu_2\Big\}.
\end{equation}
Here the supremum runs over all bounded  Lipschitz functions $\psi$ (in this case Theorem~5.10 in~\cite{villani-book2} can be extended to Lipschitz instead of continuous functions, see~\cite[Rmk.~3.6]{kuwada10}) and $Q\psi$ is the inf-convolution of $\psi,$ defined on $\X$ by
$$
Q\psi(x)=\inf_{y\in \X}\Big\{\psi(y)+\frac{d(x,y)^2 }{2}\Big\}.
$$
 The Wasserstein space $(\mathcal P_2 (\X), W_2)$ is described in the reference books~\cite{ambrosio-gigli-savare} and~\cite{villani-book2}. 
We shall define the entropy $\ent{\mu}{f}$ of a probability density with respect to a (finite or not) measure $\mu$ by $\ent{\mu}{f} = \int f \, \log f \, d\mu$ if $f \vert \log f \vert \in \mathbb{L}^1 (\mu)$ and $\infty$ otherwise.
\medskip

Our result will be stated in the two settings of a Riemannian Markov triple $(\M,\mu,\Gamma)$ ($RMT$ in short), and a Riemannian energy measure space $(\X,\tau,\mu, \mathcal{E} )$ ($REM$ in short).  These settings will be described in detail in Sections~\ref{sec-riemannian} and~\ref{sec-mms} respectively. A $REM$ space is a particular metric measure space, developed in~\cite{AGS_BE}. A $RMT$ is a smooth Riemannian manifold equipped with a weighted Laplacian (see~\cite{bgl-book}) and is a particular example of $REM$ space.

Even if a $RMT$ is a $REM$ space we prefer to state and prove our result in both settings since the argument is a little simpler in the Riemannian case. We also believe that it emphasizes the main difficulties when generalizing a result from a smooth setting  to an abstract metric measure space. In both spaces, $(P_t)_{t\geq0}$ denotes the associated Markov semigroup. It is defined through the weighted Laplacian in the $RMT$ case, and through the Dirichlet form in the $REM$ case.

The $CD(R,m)$ curvature-dimension condition is defined using the Bochner inequality~\eqref{eq-bochner} in a Riemannian manifold and in a weak form in a metric measure space (see Definitions~\ref{def-cd} and \ref{def-weak-cd}).

Recall finally that for $r\in\R$ the map $s_{r}$  is defined on $\mathbb R$ by 
$$
s_{r}(x)=
\left\{
\begin{array}{ll}
\disp \sin(\sqrt{r}\,x) / \sqrt{r} \,\,& {\rm if}\,\, r>0\\
\disp \sinh(\sqrt{|r|}\,x) / \sqrt{|r|} \,\,& {\rm if}\,\, r<0\\
x\,\,& {\rm if}\,\, r=0.
\end{array}
\right.
$$

\begin{ethm}[Equivalence between contractions and $CD(R,m)$ condition]
\label{thm-legros}
~

Consider a $RMT$ or $REM$ space as in Sections~\ref{sec-riemannian} and~\ref{sec-mms}, with (finite or not) reference measure~$\mu$ and  associated semigroup $(P_t)_{t\geq 0}$. Let $R\in\R$ and $m>0$. Then the following properties are equivalent: 
\begin{enumerate}
\item the $CD(R,m)$ (or weak $CD(R,m)$ in a $REM$ space) curvature-dimension condition holds; 
\item for any $t\geq0$ and any probability densities $f,g$ with respect to $\mu$, 
\begin{multline}
\label{eq-contraction-sh}
s_{\frac Rm}\left(\frac 12 W_2(P_tf\mu,P_tg\mu)\right)^2\leq e^{-2Rt}\,s_{\frac {R}{m}}\left(\frac {1}{2} W_2(f\mu,g\mu)\right)^2
 \\
- 2 m \int_0^t    e^{-2R(t-u)}\sinh^2 \Big( \frac{\ent{\mu}{P_uf} - \ent{\mu}{P_ug}}{2m}\Big) du;
\end{multline}
\item for  any $t\geq0$ and any probability densities $f,g$ with respect to  $\mu$,
\begin{equation}
\label{eq-contraction-square-general-2}
\!\!\!\!\!\!\!\!\!\!\!\!\!\!\!\!\!\!\!\!\!\!\!W_2^2(P_t f\mu,P_tg\mu) \! \leq \!  e^{-2Rt}W_2^2(f\mu, g\mu)
\!  -\! \frac{2}{m}\int_0^te^{-2R(t-u)} \left( \ent{\mu}{P_uf} - \ent{\mu}{P_ug} \right)^2  du.
\end{equation}
\end{enumerate} 
\end{ethm}
See Theorems~\ref{thm-main} and~\ref{thm:main_mms} 
for a more precise framework of Theorem~\ref{thm-legros}.

\smallskip

A bound with the same additional term as in (ii) has also been derived in~\cite{BGG15} for some specific instances of symmetric Fokker-Planck equations in $\R^m$, for which the generator only satisfies a $CD(R,\infty)$ condition. Combined with a deficit in the Talagrand inequality, it has led to refined convergence estimates on the solutions.

Next section presents the easiest part of the proof of Theorem~\ref{thm-legros}. More precisely, $(\rm i)\Rightarrow (ii)$ and an outline of $(\rm ii) \Rightarrow (iii)$, including the key~Proposition~\ref{prop-3-4}. The full proof of $(\rm ii) \Rightarrow (iii)$ requires some knowledge on the spaces and will be finished in Section~\ref{subsec:3-4}. Sections~\ref{sec-proofdebut} to~\ref{sec-mms} (but Section~\ref{subsec:3-4}) are dedicated to the more difficult $(\rm iii)\Rightarrow (i),$ in both $RMT$ and $REM$ spaces.

\section{Proof of Theorem~\ref{thm-legros}: first implications}
\label{sec-proofdebut}

%{\red We first present some implications of Theorem~\ref{thm-legros}, which are easy or follow from earlier results.}
\subsubsection*{Proof of $(\rm i)\Rightarrow (ii)$} 

In~\cite{EKS13} M. Erbar, K.-T. Sturm and the fourth author of this paper have proved an {\it Evolutional variational inequality} (EVI in short) in the $REM$ spaces.  Let $f,g$ be probability densities with respect to  $\mu$ and let $U_m = \exp ( -  \ent{\mu}{\cdot}/m )$. Then, under the weak $CD(R,m)$ condition, 
\begin{equation}\label{eq-evi-eks}
\frac{d}{dt} s_{\frac Rm}\left(\frac12 W_2(P_tf\mu,g\mu) \right)^2 + R \, s_{\frac Rm}\left(\frac12 W_2(P_tf\mu, g\mu) \right)^2 \leq \frac{m}{2} \left( 1 - \frac{U_m(g)}{U_m(P_t f)} \right).
\end{equation}
But it is classical, see e.g.~\cite{ambrosio-gigli-savare}, how to deduce a contraction property in $W_2$ distance between solutions $(P_tf)_{t \geq 0}$ and $(P_tg)_{t \geq 0}$ from an EVI: one applies the EVI to the curve  $(P_tf)_{t \geq 0}$ and $P_sg$ for a given $s$, and then (with the time variable $s$) to the curve $(P_sg)_{s \geq 0}$ and $P_tf$ for a given $t$; then one adds both inequalities,
takes $t=s$ and integrate in time. Then one obtains (ii). To sum up, it turns out that the EVI~\eqref{eq-evi-eks} not only leads to the property~\eqref{eq-eks}, as observed in~\cite{EKS13}, but also to the {\it same-time} contraction property~(ii).

\subsubsection*{Outline of the proof of $(\rm ii)\Rightarrow (iii)$}

We first observe that $\sinh^2(x)\geq x^2$ for any $x$, so ii) implies the same bound with $\sinh^2(x)$ replaced by $x^2$ in the integral. Then the
 implication $(\rm ii)\Rightarrow (iii)$ is essentially a consequence of the following result, which we prove in the general context of a geodesic space. 

\begin{eprop}
\label{prop-3-4}
Let $( Y, d_Y )$ be a geodesic metric space, $U : Y \to ( - \infty , \infty ]$ and 
$\varphi_t : Y \to Y$ $(t \ge 0)$ a one-parameter family of maps. 
Suppose that $t \mapsto \varphi_t (y)$ is continuous for all $y \in Y$ and 
$U ( \varphi_t (y) ) \in \R$ for all $t > 0$ and $y \in Y.$
Suppose also that for $y_0 , y_1 \in Y$ and $t > 0$, 
\begin{multline}
\label{contrtt2}
s_{\frac Rm} \left( 
    \frac12 d_Y (\varphi_t (y_0) , \varphi_t (y_1) )
\right)^2
\leq 
e^{-2Rt}\,
s_{\frac {R}{m}} \left(
    \frac {1}{2} d_Y ( y_0 , y_1 )
\right)^2 
\\
- \frac{1}{2 m} \int_0^t 
e^{-2R(t-u)} 
\PAR{ 
    U ( \varphi_u (y_0) ) - U ( \varphi_u (y_1) )
} ^2
du.
\end{multline} 
Then 
\begin{equation*}
d_Y (\varphi_t (y_0) , \varphi_t (y_1) )^2
\leq 
e^{-2Rt}\, d_Y ( y_0 , y_1 )^2
\\
- \frac{2}{m} \int_0^t 
e^{-2R(t-u)} 
\PAR{ 
    U ( \varphi_u (y_0) ) - U ( \varphi_u (y_1) )
} ^2
du.
\end{equation*}
\end{eprop}

\begin{Proof}
We adapt the argument of \cite[Prop.~2.22]{EKS13}.  Let $( y_s )_{s \in [ 0 , 1 ]}$ be a geodesic from $y_0$ to $y_1$ in $Y$, and let $t>0$ be fixed. For any $n$ and $1 \leq i \leq n, $ let $x_i^n = d_Y ( \varphi_t ( y_{(i-1)/n} ), \varphi_t ( y_{i/n} ) )$. Then
$$
d_Y ( \varphi_t (y_0) , \varphi_t (y_1) )^2
\leq 
\Big( \sum_{i=1}^n x_ i^n \Big)^2
\leq
n \sum_{i=1}^n (x_i ^n)^2
$$
for any $n.$ In particular
$$
d_Y ( \varphi_t (y_0) , \varphi_t (y_1) )^2
\leq \limsup_{n \to \infty} n \sum_{i=1}^n (x_i ^n)^2.
$$
Now, by neglecting the second term in the right-hand side of~\eqref{contrtt2} and by geodesic property,
$$
s_{\frac{R}{m}} \Big(\frac{x_i ^n}{2} \Big) \leq e^{-Rt} \, s_{\frac{R}{m}} \Big( \frac{1}{2} d_Y (y_{(i-1)/n} , y_{i/n}) \Big)  = e^{-Rt} \, s_{\frac{R}{m}} \Big( \frac{1}{2n} d_Y ( y_0, y_1) \Big).
$$
It follows, as in \cite[(2.32)]{EKS13}, that there exists a constant $c$ such that $x_i^n \leq c/n$ for large $n$ and any $1 \leq i \leq n.$ Moreover  $ s_{\frac{R}{m}} (x)^2 = x^2- R x^4/(3m)+O(x^6)$ as $x$ tends to $0$, so that
\begin{equation}\label{lemmaDLs}
\limsup_{n \to \infty}n\sum_{i=1}^n (x_i^n)^2= 4 \limsup_{n \to \infty}n\sum_{i=1}^n s_{\frac{R}{m}}(x_i^n/2)^2.
\end{equation}
As a consequence
\begin{align*}
d_Y ( \varphi_t ( y_0 ) , \varphi_t ( y_1 ) )^2
& \leq 
4 \limsup_{n \to \infty} 
n \sum_{i=1}^{n} s_{\frac{R}{m}} 
\left( 
    \frac12 d_Y ( \varphi_t ( y_{(i-1)/n} ) , \varphi_t ( y_{i/n} ) )
\right)^2 
\\
& \leq 
4 \limsup_{n \to \infty} \Bigg(
n \sum_{i=1}^n e^{-2Rt}\,
s_{\frac {R}{m}} \left(
    \frac {1}{2} d_Y ( y_{(i-1)/n} , y_{i/n} )
\right)^2 
\\
& \hspace{4em}
- \frac{1}{2m} \int_0^t 
e^{-2R(t-u)} 
n \sum_{i=1}^n  
\PAR{ 
    U ( \varphi_u ( y_{(i-1)/n} ) ) - U ( \varphi_u ( y_{i/n} ) ) 
}^2 
du
\Bigg)
\end{align*}
by assumption~\eqref{contrtt2}. 

Then the conclusion follows from this estimate 
by using~\eqref{lemmaDLs} with $d_Y ( y_{(i-1)/n} , y_{i/n})$ in place of $x_i^n$ in the first term, and the Cauchy-Schwarz inequality in the second term. 
\end{Proof}

 Let us return to our case. 
As stated before, we give an outline and 
leave the rigorous argument to Section~\ref{subsec:3-4}. 
We apply Proposition~\ref{prop-3-4} for $( Y , d_Y ) = ( \mathcal{P}_2 (\X) , W_2 )$. 
Under (ii), we can extend the action of the semigroup $P_t$ to probability measures. 
Then, $\varphi_t = P_t$ fulfills all the assumptions of Proposition~\ref{prop-3-4} 
with $U = \entf{\mu}$. 
This ensures that (ii) implies~(iii). 
\medskip

To sum up, (i) implies (ii), and (ii) implies (iii) 
(with the aid of an additional argument in Section~\ref{subsec:3-4}). 
Thus, it remains to prove  that, conversely, (iii) implies (i).

%%%%%%%%%%%%%%%%%%%%
\section{Strategy of the converse proof}
\label{sec-strategy}
The proof of (iii) $\Rightarrow$ (i) will be given in the two cases of a $RMT$ (in Section~\ref{sec-riemannian}) and a $REM$ space (in Section~\ref{sec-mms}).  In this section we present its strategy, in a formal way.

\subsection{Example of a gradient flow in $\R^d$}
\label{sec-radient-flow}

Let us first present the easiest case of a smooth gradient flow in $\R^d$.  There we shall see that the equivalence between the contraction inequality~\eqref{eq-contraction-square-general-2} and the $CD(R,m)$ curvature-dimension condition is natural. It gives a way to understand the general case. 

\medskip

Let $F : \R^d \to \R$ be a $\mathcal C^2$ smooth function, and let $(X_t)_{t\geq0}$ be a gradient flow for the function $F$, that is, a solution to the differential equation 
\begin{equation}
\label{eq-fg}
\frac{dX_t}{dt}=-\nabla F(X_t). 
\end{equation}

Following~\cite{EKS13}, the function $F$ satisfies a $CD(R,m)$ curvature-dimension condition for $R\in\R$ and $m>0$ if for any $x,h\in\R^d$, the map  $[0,1]\ni s\mapsto \phi(s)=F(x+sh)$ satisfies the convexity inequality
\begin{equation}
\label{eq-cd-R}
\phi''(s)\geq R || h ||^2 +\frac{1}{m}(\phi'(s))^2.
\end{equation}
Here $||\cdot||$ is the Euclidean norm in $\R^d$. Since  the path $(x+sh)_{s\in[0,1]}$ is a geodesic between $x$ and $x+h$, this means that
$F$ satisfies a $(R,m)$-convexity condition along geodesics. 
 
Let now $(X_t)_{t\geq0}$ and $(Y_t)_{t\geq0}$ be two solutions to~\eqref{eq-fg} with initial conditions $X_0$ and $Y_0$ respectively. The function
 $\Lambda(t)=||X_t-Y_t||^2$ satisfies
 $$
\Lambda'(u)=-2\int_0^1\phi_u''(s)ds 
$$ 
where $\phi_u(s)=F(X_u+s(Y_u-X_u))$. 
If now the function $F$ satisfies the above $CD(R,m)$ condition~\eqref{eq-cd-R}, then 
$$
\Lambda'(u)
\leq 
-2R||X_u-Y_u||^2
-\frac{2}{m}\int_0^1 (\phi_u'(s))^2 du
\leq 
-2R\Lambda(u)-\frac{2}{m} (\phi_u(1)-\phi_u(0))^2. 
$$ 
Integrating over the interval $[0,t]$, we get 
\begin{equation}
\label{eq-contraction-R}
||X_t-Y_t||^2\leq e^{-2Rt}||X_0-Y_0||^2-\frac{2}{m}\int _0^te^{-2R(t-u)}(F(X_u)-F(Y_u))^2du. 
\end{equation}

Conversely, let us assume that the gradient flow driven by $F$ satisfies the property~\eqref{eq-contraction-R} for any $t\geq0$ and any initial conditions $X_0$ and $Y_0.$ Then  $F$ satisfies the $CD(R,m)$ condition~\eqref{eq-cd-R}. For, taking the time derivative of~\eqref{eq-contraction-R} at $t=0$ implies 
$$
- (X_0-Y_0)\cdot(\nabla F(X_0)-\nabla F(Y_0))\leq -R||X_0-Y_0||^2-\frac{1}{m}(F(X_0)-F(Y_0))^2. 
$$
Let then $x, h$ in $\R^d$ and $s \in [0,1]$ be fixed. A Taylor expansion for $Y_0 = x+(s+\varepsilon) h$ tending to $X_0 = x + sh$ (along a geodesic), so for $\varepsilon \to 0,$  implies back the  $CD(R,m)$ condition~\eqref{eq-cd-R}.  

\medskip

Let us observe that inequality~\eqref{eq-contraction-R} is exactly~\eqref{eq-contraction-square-general-2} when replacing $\R^d$ with the space of probability densities, the Euclidean norm with the Wasserstein distance, $F$ with the entropy, $(X_t)_{t\geq0}$ with the semigroup $(P_t)_{t \geq 0}$ and the $CD(R,m)$ condition~\eqref{eq-cd-R} with the corresponding Bakry-\'Emery condition, which is equivalent to the $(R,m)$-convexity of the entropy (see \cite{EKS13}).
Of course, this computation is natural since
% at considered since~\cite{jko}
the considered evolution is the gradient flow of the entropy with respect to the Wasserstein distance, see~\cite{ambrosio-gigli-savare,jko}.
\smallskip

We now want to mimic the above proof for a smooth gradient flow on $\R^d$ to the setting of a general semigroup on $(\mathcal P_2 (\X), W_2)$.
As here in the smooth case, we shall see in the coming section that geodesics play a fundamental role.

%%%%%%%%%%%%%%%%%%%%
\subsection{How to adapt the gradient flow proof to the general case?}
\label{sec-how}

The most natural method to prove that a contraction inequality in Wasserstein distance, as in~\eqref{eq-premiere}, implies a curvature condition is to use close Dirac measures as initial data (see e.g.~\cite{bgl-15}). In our case, this can not be performed since the entropy of a Dirac measure is infinite. There seems to be hope since we consider the entropy of the heat kernel 
in positive time, when it becomes finite. 
However, it does not work again if we are on a homogeneous space. 
For instance, on $\mathbb{R}^d$, the entropy of the heat kernel 
$p_t (x,\cdot)$ does not depend on $x$ and the dimensional corrective terms 
in Theorem~\ref{thm-legros} vanish 
if we consider two Dirac measures as initial data.

 To solve this issue we shall consider as initial data a probability density $g$ (with respect to $\mu$) and a perturbation of it, both in sufficiently wide classes of functions. The perturbation will be built by means of a geodesic in the Wasserstein space $(\mathcal P_2 (\X), W_2)$. More precisely, given such a $g$,
we are looking for  a path $(g_s)_{s \geq0}$ of probability densities whose Taylor expansion for small $s$ is a geodesic in $\mathcal P_2 (\X)$ with a direction given by a function~$f$. 
We explain the idea on a $RMT$.

For that, consider the generator $L^g=L + \Gamma(\log g,\cdot)$ (see~\eqref{eq-Gamma} for the definition of $\Gamma$) with associated semigroup $(P_t^g)_{t\geq0}$. Given a direction function $f$, there are two ways of defining the path $(g_s)_{s \geq 0}$, both admitting the same Taylor expansion for small $s$:
\begin{itemize}
\item One can first consider the path  $g_s=g(1-sL^g f)$ for small $s$ and a smooth and compactly supported function $f$.  The function $g_s$ is a smooth, bounded and compactly supported perturbation of $g$. This path will be used on a $RMT$ since such functions are adapted to the Riemannian setting. 
\item One can also consider the path $\tilde{g}_s=g(1+f-P_s^gf)$, again for $s$ small and ``nice'' $f\in\mathbb L^\infty(\mu)$. The path $(\tilde{g}_s)$ has the same Taylor expansion as $(g_s)$ since  $f-P_s^gf=-sL^gf+o(s)$. 
This path will be used on $REM$ spaces. Indeed, regularity of functions (such as $g_s$ above) is clearly a difficult issue in the setting of metric measure spaces, and $\mathbb L^\infty(\mu)$ functions are much more adapted to them. By using the semigroup $( P_s^g )_{s \ge 0}$ instead of the generator $L^g$, we can apply the maximum principle which preserves (essential) boundedness of functions.
\end{itemize}

\begin{erem} \label{rem:Approximation_geodesic}
Let us see, formally and in the Euclidean space $\R^d$, why the probability measure $ g_s dx$ has the same first-order Taylor expansion as the geodesic in the Wasserstein space. Let $\nu_0$ be a probability measure in $\R^d$ being absolutely continuous with respect to the Lebesgue measure, $\psi:\R^d\rightarrow\R$ be a convex map, and
$$
\nu_s=((1-s){\rm Id}+s\nabla\psi)_{\#}\nu_0
$$
for $s\in[0,1]$. The path $(\nu_s)_{s\in[0,1]}$ is a geodesic path between $\nu_0$ and $\nu_1$ in the Wasserstein space, that is for any $s,t\in[0,1]$, 
$$
W_2(\nu_s,\nu_t)=|t-s| \, W_2(\nu_0,\nu_1).
$$
Moreover, for any test function $H:\R^d\rightarrow \R$, and by a formal Taylor expansion when $s$ goes to~0,
$$
\int Hd\nu_s=\int H((1-s)x+s\nabla\psi(x))d\nu_0(x)=\int\SBRA{ H(x)+s\nabla H(x)\cdot (\nabla\psi(x)-x)+o(s)}d\nu_0(x). 
$$
Assume now that $d\nu_0=gdx$ for a function $g$. Then, by integration by parts as in~\eqref{eq-ipp} below,
$$
\int Hd\nu_s=\int Hd\nu_0-s\int H \, L^g(f) \, d\nu_0+o(s) = \int H \, g_s dx + o(s)
$$
where $f(x)=\psi(x)-|x|^2/2$. 

In conclusion, the  path $(g_s)_{s \geq 0}$ appears as a (smooth) first-order Taylor expansion of the $W_2$-geodesic path  $(\nu_s)_{s \geq 0}$.
\end{erem}

%%%%%%%%%%%%%%%%%%%%
\section{The Riemannian Markov triple context}
\label{sec-riemannian}

In this section we prove the implication (iii) $\Rightarrow$ (i) of Theorem~\ref{thm-legros} in the context of a Riemannian manifold, in the form of Theorem~\ref{thm-main} below.  

\subsection{Framework and results}

Let $(\M, \mathcal G)$ be a  connected complete $\mathcal C^\infty$-Riemannian manifold. Let $V$ be a $\mathcal C^{\infty}$ function on $\M$ and consider the Markov semigroup $(P_t)_{t\geq 0}$ with generator $L=\Delta-\nabla V\cdot\nabla$,  where $\Delta$ is the Laplace-Beltrami operator.  Let also $ d\mu=e^{-V}dx$ where $dx$ is the Riemannian measure and 
$\Gamma$ be the carr\'e du champ operator, defined by
\begin{equation} \label{eq-Gamma}
\Gamma(f,g)=\frac{1}{2}(L(fg)-fLg-gLf) 
\end{equation}
for any smooth $f,g.$ We let $\Gamma(f)=\Gamma(f,f)=|\nabla f|^2$ where  $|\nabla f|$ stands for the length of $\nabla f$ with respect to the Riemannian metric $\mathcal G$. 

\medskip

Then $(\M,\mu,\Gamma)$ is a full Markov triple in a Riemannian manifold, as in~\cite[Chap.~3]{bgl-book}, and  in this work we call it a {\it Riemannian Markov triple $(RMT)$}.

\medskip

The measure $\mu$ is reversible with respect to the semigroup,
that is, for any $t\geq0$, $P_t$ is a self-adjoint operator in $\dL^2(\mu)$.  
Moreover the integration by parts formula
$$
\int fLg \, d\mu=-\int \Gamma(f,g)d\mu 
$$
holds for all $f, g$ in the set $\mathcal{C}^\infty_c (\M)$ of infinitely differentiable and compactly supported functions on~$\M$.
The generator $L$ satisfies the diffusion property, that is, for any smooth functions $\phi, f,g$, 
$$
L(\phi(f))=\phi'(f)Lf+\phi''(f)\Gamma(f), 
$$
or equivalently 
\begin{equation} \label{eq:diffusion2}
\Gamma(\phi(f),g)=\phi'(f)\Gamma(f,g). 
\end{equation}
In other words, the carr\'e du champ operator is a derivation operator for each component. 
\medskip

The map $(x,t)\mapsto P_tf(x)$ is simply the solution to the parabolic equation $\partial_t u=Lu$ with $f$ as the initial condition. 

\begin{edefi}[$CD(R,m)$ condition]
\label{def-cd}
Let $R\in\R$ and $m\in (0,\infty]$. We say that the $RMT$ $(\M,\mu,\Gamma)$  satisfies a  $CD(R,m)$ curvature-dimension condition if 
 $$
\Gamma_2(f)\geq R\Gamma(f)+\frac{1}{m}(Lf)^2
$$
for any smooth function $f,$ say in $\mathcal C^\infty_c(\M),$
where
\begin{equation} \label{eq:Gamma2}
\Gamma_2(f)=\frac{1}{2}(L\Gamma(f)-2\Gamma(f,Lf)).
\end{equation}
\end{edefi}

Let us notice that $m$ can be different from the dimension of the manifold $\M$.
The $CD(R,m)$ curvature-dimension condition is called the Bakry-\'Emery or $\Gamma_2$ condition 
and has been introduced in~\cite{bakryemery}
 (see also~the recent~\cite{bgl-book}).

\begin{eex}
On a $d$-dimensional Riemannian manifold $(\M, \mathcal G)$
\begin{itemize}
\item the operator $L=\Delta$ satisfies a $CD(R,m)$ condition if $m\geq d$ and the Ricci curvature of the manifold is bounded from below by $R$; 

\item more generally, the operator  $L=\Delta-\nabla V\cdot\nabla$ satisfies a $CD(R,m)$ condition if  $m \geq d$ and
$$
\mathrm{Ric}+{\rm Hess}(V)\geq R \, \mathcal G+\frac{1}{m-d}\nabla V\otimes\nabla V, 
$$
where $\mathrm{Ric}$ is the Ricci tensor of $(\M, \mathcal G)$, see for instance~\cite[Sec.~C6]{bgl-book} (when $m=d$ then we need $V=0$).  
\end{itemize}
\end{eex}

In a $RMT$, the following result gives the implication (iii) $\Rightarrow$ (i) in Theorem~\ref{thm-legros} :

\begin{ethm}
\label{thm-main} Let $(\M,\mu,\Gamma)$ be a Riemannian Markov triple and $(P_t)_{t\geq0}$ its   associated Markov semigroup. Let $R\in\R$ and $m>0$.  
If the inequality~\eqref{eq-contraction-square-general-2} holds for any $t\geq0$ and  any smooth functions $f,g$ on $\M$ with $f \mu, g\mu$ in $\mathcal P_2(\M)$, 
then the $CD(R,m)$ condition of Definition~\ref{def-cd} holds.
\end{ethm}

\subsection{Proof of Theorem~\ref{thm-main}}
\label{sec-proof-riemannian}

It is based on the approximation of geodesics introduced in Section~\ref{sec-how} (see Remark~\ref{rem:Approximation_geodesic}), properties of the Hopf-Lax solution of the Hamilton-Jacobi equation, and an adapted class of test functions.

\medskip

Let $f$ be in $\mathcal C^\infty_c(\M).$ 
Let also $g$ be a smooth and positive function on $\M$ such that $g\mu\in\mathcal P_2(\M),$
$$
\int g \, |\log g| \, d\mu<\infty\quad{\rm and}\quad\int \frac{\Gamma(g)}{g}d\mu<\infty.
$$ 
Let us define the generator $L^g$ by
$$
L^g h =Lh+\Gamma(\log g,h)
$$
on smooth functions $h$. Since $g>0$, then $L^g$ is well defined on the set $\mathcal C^\infty_c(\M)$ and $L^g h \in \mathcal C_c^\infty(\M)$ for any $h\in\mathcal C^\infty_c(\M)$. Moreover, the  generator $L^g$ satisfies an integration by parts formula with respect to the probability measure $g\mu$ : for $h,k \in \mathcal{C}^\infty_c (\M)$ (one of them can be with non compact support)
\begin{equation}
\label{eq-ipp}
\int h \, L^g k \, gd\mu=-\int \Gamma(h,k) \, gd\mu.
\end{equation}

For any $s\geq 0$, let us define $g_s=g(1- sL^gf)$. The function $L^g f$ is in $\mathcal C^\infty_c(\M)$, so bounded, and we can let $N=||L^gf||_\infty.$ We shall frequently use the bounds $(1-sN) g \leq g_s \leq (1+sN)g$. In particular $g_s>0$ for $s< 1/N.$ Moreover $\int g_sd\mu=1$. Hence, for $s$ small enough, which we now assume, $g_s\mu$ is in $\mathcal P_2(\M)$ with a smooth and positive density.
The proof of Theorem~\ref{thm-main} consists in applying~\eqref{eq-contraction-square-general-2} with $g_s$  instead of $f,$ dividing by $2s^2$ and letting $s$ go to $0$. For this we shall estimate the three terms in the inequality.

\medskip

 A key tool is the Hopf-Lax semigroup defined on bounded Lipschitz functions $\psi$ by
\begin{equation} \label{eq:Hopf-Lax}
Q_s\psi(x):=\inf_{y\in \M}\BRA{\psi(y)+\frac{d(x,y)^2}{2s}}, \quad s>0, \,\,x \in \M.
\end{equation}
The map  $x\mapsto Q_s\psi(x)$ is Lipschitz for every $s\geq0$, and the map $(s,x)\mapsto Q_s\psi(x)$
satisfies  the Hamilton-Jacobi equation
$$
\partial_s Q_s\psi+\frac{1}{2}|\nabla Q_s\psi|^2=0,\quad \lim_{s\rightarrow 0}Q_s\psi=\psi
$$
in a sense given in~\cite[Thms. 22.46 and~30.30]{villani-book2} for instance.
We observe that $s Q_s(\psi) = Q_1(s\psi) = Q(s\psi)$, so for $s>0$ the Kantorovich duality~\eqref{eq-kanto} can be written as 
\begin{equation}
\label{eq-23}
\frac{W_2^2(\nu_1, \nu_2)}{2s^2} = \frac{1}{s} \sup_{\psi} \SBRA{\int Q_s\psi \, d\nu_1 -\int \psi \, d\nu_2}.
\end{equation}

\medskip

\noindent{\bf Estimate on the term on the left-hand side of~\eqref{eq-contraction-square-general-2}.}  
Letting $\psi=f$ in~\eqref{eq-23}, we obtain  
\begin{equation} \label{eq:first-estimate0}
\frac{W_2^2(P_tg_s\mu,P_tg\mu)}{2s^2}\geq \int \frac{Q_s f P_tg_s- f P_tg}{s} d\mu. 
\end{equation} 
Since $f$ is Lipschitz, almost everywhere in $\M$ we have  
$$
\lim_{s\rightarrow0}\frac{Q_s f P_tg_s- f P_tg}{s} = -\frac{1}{2}\Gamma(f)P_tg-fP_t(gL^gf)
$$
by (vii') in~\cite[Thm.~30.30]{villani-book2}.  But, by the definition of $Q_s f$ and since $f$ is bounded, 
\[
Q_s f (x) 
= 
\inf_{ y \in B ( x, \sqrt{ 4s \| f \|_\infty } )} 
\left\{ 
f (y) + \frac{ d (x,y)^2 }{ 2s } 
\right\}. 
\]
Thus, for the Lipschitz seminorm $\Vert \cdot \Vert_{Lip}$,
\begin{eqnarray}\label{minoqsf}
0 \geq \frac{Q_s f(x)  - f (x)}{s}
& \geq 
\displaystyle \inf_{ y \in B ( x, \sqrt{ 4s \| f \|_\infty } ) \setminus \{ x \} } 
\left\{ 
\frac{ f (y) - f (x) }{ d(x,y) } \frac{d (x,y)}{s} + \frac{ d (x,y)^2 }{ 2s^2 } 
\right\} \nonumber
\\
& \geq
\displaystyle - \frac{1}{2} \sup_{ y \in B ( x, \sqrt{ 4s \| f \|_\infty } ) \setminus \{ x \} } 
\left( 
\frac{ f (y) - f (x) }{ d(x,y) } 
\right)^2 
\geq 
- \frac12 \| f \|_{Lip}^2  
\end{eqnarray}
(see also~\cite[page 585]{villani-book2}). Moreover $|| Q_s f||_\infty \leq ||f||_\infty,$ so, adding and subtracting $Q_s f P_t g$,
$$
\Big|\frac{Q_s f P_tg_s- f P_tg}{s}\Big|\leq || Q_s f||_\infty \; \vert P_t(g L^gf) \vert +P_tg \; \frac{f-Q_sf}{s}\leq \Big( ||f||_\infty \, ||L^g f||_\infty +\frac{||f||^2_{Lip}}{2} \Big) P_tg.
$$
The right-hand side is in $\mathbb L^1(\mu)$, so by the Lebesgue dominated convergence theorem   
$$
\liminf_{s\rightarrow0}\frac{W_2^2(P_tg_s\mu,P_tg\mu)}{2s^2}\geq \int \PAR{-\frac{1}{2}\Gamma(f)P_tg-fP_t(gL^gf)}d\mu.
$$
 Now, by reversibility of the measure $\mu$ and the  integration by parts formula~\eqref{eq-ipp},
$$
\int fP_t(gL^gf)d\mu=\int P_tf L^g(f) \, gd\mu=-\int\Gamma(f, P_tf) \, gd\mu.
$$ 
Thus we obtain our first estimate:
\begin{equation}
\label{eq-first-estimate}
\liminf_{s\rightarrow 0} \frac{W_2^2(P_tg_s\mu,P_tg\mu)}{2s^2}\geq -\frac{1}{2}\int P_t(\Gamma(f)) gd\mu+\int \Gamma(f,P_tf)gd\mu.
\end{equation}
{\bf Estimate on the first term on the right-hand side}. According to~\eqref{eq-23} we need an upper bound on the quantities $\int Q_s(\psi) g_sd\mu-\int \psi gd\mu$, independent of the bounded Lipschitz function~$\psi$. 

\medskip

 First of all, for $0 < t < s$,
 \begin{equation} \label{eq-hl}
\frac{d}{dt}\int Q_t \psi \, g_t \, d\mu=\int \SBRA{-\frac{1}{2}\Gamma(Q_t\psi)(1-tL^gf)-Q_t\psi L^g f}gd\mu. 
\end{equation}
This is justified by item (vii) in~\cite[Thms. 22.46 and~30.30]{villani-book2} and the properties that $g\mu\in\mathcal P(\M),$ $L^g f$ is bounded, $||Q_t \psi ||_{\infty}\leq ||\psi||_{\infty}$ and $||Q_t \psi ||_{Lip}\leq ||\psi||_{Lip}$ for any $t$. 

 Now the integration by parts formula~\eqref{eq-ipp} gives $-\int Q_t\psi \,L^g f \,gd\mu=\int \Gamma(Q_t\psi,f)gd\mu$. Recall that $L^gf$ is bounded 
 and that we have let $N=||L^gf||_\infty$.  For $t<s< 1/N$ we obtain
\begin{multline*}
\frac{d}{dt}  \int Q_t \psi \, g_t \, d\mu \leq \int \SBRA{-\frac{1}{2}\Gamma(Q_t\psi)(1-sN)+\Gamma(Q_t\psi,f)}gd\mu
\\
=
\int \SBRA{-\frac{1-sN}{2}\Gamma\PAR{Q_t\psi-\frac{1}{1-sN}f}+\frac{1}{2(1-sN)}\Gamma(f)}gd\mu\leq
  \frac{1}{2(1-sN)}\int\Gamma(f)gd\mu. 
\end{multline*}
Integrating over the set $t\in[0,s]$ :
$$
\int Q_s \psi \, g_sd\mu-\int \psi gd\mu \leq  \frac{s}{2(1-sN)}\int\Gamma(f)gd\mu.
$$

Finally the Kantorovich duality~\eqref{eq-23} gives our second estimate: 
\begin{equation}
\label{eq-second-estimate}
\limsup_{s\rightarrow 0}\frac{W_2^2(g_s\mu,g\mu)}{2s^2}\leq\frac{1}{2}\int\Gamma(f)gd\mu. 
\end{equation}
{\bf Estimate on the second term on the right-hand side}. Let $u > 0$ and let us compute the limit of $\frac{1}{s}\PAR{\ent{\mu}{P_ug_s}-\ent{\mu}{P_ug}}$ when $s$ goes to 0.
First, for any $s > 0$, 
$$
\frac{d}{ds}P_u(g_s)\log P_u(g_s)=-(1+\log P_ug_s) \; P_u \big( g L^g  f \big).
$$
Then, for $0 < s < 1/N$, 
$$
|(1+\log P_ug_s)P_u(g L^g f)|\leq NP_ug \; (1 + \log(1+N) + |\log P_u(g)|). 
$$
Forgetting the dimensional corrective term in~\eqref{eq-contraction-square-general-2}, by the von Renesse-Sturm theorem~\cite{sturm-vonrenesse} the $RMT$ satisfies a $CD(R,\infty)$ condition. In particular, and since $\int \Gamma(g) / g \, d\mu<\infty,$ one can use a local logarithmic  Sobolev inequality~\cite[Thm.~5.5.2]{bgl-book} to deduce
$\int P_ug \, |\log P_u g| \, d\mu<\infty$.  In particular the right-hand side in the last inequality is in $L^1(\mu).$ 
Then, by the Lebesgue convergence theorem and~\eqref{eq-ipp},
\begin{eqnarray*}
\lim_{s\rightarrow0}\frac{{\ent{\mu}{P_ug_s}-\ent{\mu}{P_ug}}}{s}
&=&
-\int(1+\log P_ug)P_u \big(g L^g  f  \big)d\mu
\\
&=&
-\int P_u(\log P_ug)L^g  f \; gd\mu =\int \Gamma(P_u(\log P_ug),f)gd\mu. 
\end{eqnarray*}
By the Fatou lemma we obtain the third estimate : 
\begin{multline}
\label{eq-third-estimate}
\limsup_{s\rightarrow0}-\frac{1}{m}\int_0^te^{-2R(t-u)}\SBRA{\frac{\ent{\mu}{P_ug_s}-\ent{\mu}{P_ug}}{s}}^2du
\\
\leq
-\frac{1}{m}\int_0^te^{-2R(t-u)}\PAR{\int \Gamma(P_u(\log P_ug),f)gd\mu}^2du . 
\end{multline}
{\bf Conclusion.} Dividing the inequality~\eqref{eq-contraction-square-general-2} by $2s^2$, letting $s$ go to $0$ and using the three estimates~\eqref{eq-first-estimate},~\eqref{eq-second-estimate} and~\eqref{eq-third-estimate} we get
\begin{multline*} 
-\frac{1}{2}\int P_t \Gamma(f) \, gd\mu+\int \Gamma(f,P_tf)gd\mu
\\
\leq
\frac{e^{-2Rt}}{2}\int\Gamma(f)gd\mu-\frac{1}{m}\int_0^te^{-2R(t-u)}\PAR{\int\Gamma(P_u(\log P_ug),f)gd\mu}^2du .
\end{multline*}
This inequality is an equality when $t=0$, and since $f\in\mathcal C_c^\infty(\M)$, its derivative at $t=0$ implies 
$$
-\frac{1}{2}\int L \Gamma(f) \, gd\mu+\int \Gamma(f,Lf)gd\mu\leq -R\int\Gamma(f)gd\mu-\frac{1}{m}\PAR{\int\Gamma(\log g,f)gd\mu}^2. 
$$
Since $\int\Gamma(\log g ,f)gd\mu=\int\Gamma(g,f)d\mu=-\int gLfd\mu$ and by definition of the $\Gamma_2$ operator we get
\begin{equation} \label{eq:preCD}
\int \Gamma_2(f)gd\mu\geq  R\int\Gamma(f) \, gd\mu+\frac{1}{m}\PAR{\int Lf \, g d\mu}^2
\end{equation}
for any $f\in\mathcal C_c^\infty(\M)$ and any positive smooth probability density $g$ with $\ent{\mu}{g}\!, \int \Gamma(g)/g <\infty$. 

\medskip

Inequality~\eqref{eq:preCD} appears as a weak form of the $CD(R,m)$ condition.  Again from the $CD(R,\infty)$ condition, it is a consequence of Wang's Harnack inequality (see~\cite[Thm.~5.6.1]{bgl-book} and~\cite{wang-book}) that there exist $\alpha_0 >0$ and $o \in\M$ such that 
\begin{equation} \label{eq:integrable}
 \int \exp ( - \alpha_0 d ( o, x )^2 ) \, d \mu (x) < \infty. 
\end{equation} 
Then, in~\eqref{eq:preCD} we can replace $g$ by a sequence $(g_p)_p$ converging to 
the Dirac measure $\delta_x$  at $x\in\M$; we get
$$
\Gamma_2(f)\geq R\Gamma(f)+\frac{1}{m}(Lf)^2
$$
at any $x\in\M$ and for any function $f\in\mathcal C_c^\infty(\M)$. This is  the $CD(R,m)$ condition as in Definition~\ref{def-cd}, and this finishes the proof of Theorem~\ref{thm-main}. 

%%%%%%%%%%%%%%%%%%%%%%%%%%%%%%%%%
\section{The Riemannian energy measure space context} 
\label{sec-mms}

In this section we prove the implication (iii) $\Rightarrow$ (i) of Theorem~\ref{thm-legros} in the context of a Riemannian energy measure ($REM$) space, a particular case of metric measure spaces, see Theorem~\ref{thm:main_mms} below. 
The proof goes along the same overall strategy as in the manifold case of Section~\ref{sec-proof-riemannian}. However, to overcome the lack of differentiability, it will require several tools and results from optimal transport and heat distributions on metric measure spaces.

The framework and the main Theorem~\ref{thm:main_mms} are stated in Section~\ref{subsec:frame_mms}. As an intermezzo, in Section~\ref{subsec:3-4} we complement the proof of (ii) $\Rightarrow$ (iii) in Theorem~\ref{thm-legros}.
The path $(\tilde{g}_s)_{s\geq 0}$ is constructed in Section~\ref{subsec-path}, the three key estimates are given in Section~\ref{sec-3-estimates}, finally the main proof is given in Section~\ref{subsec:pf_mms}.

\subsection{Framework and results}
\label{subsec:frame_mms}

As a natural framework, we will state our result 
on a Riemannian energy measure space, as introduced in~\cite{AGS_BE}. 
Let $( \X , \tau )$ be a Polish topological space 
and $\mu$ a locally finite Borel measure with a full support.
Let $( \mathcal{E} , \mathcal{D} (\mathcal{E} ) )$ be 
a strongly local symmetric Dirichlet form on $\mathbb{L}^2 ( \mu ).$ Let finally
$(P_t)_{t \geq 0}$ be its associated semigroup and $L$ its generator, with domain $\mathcal D(L)\subset\mathbb{L}^2 ( \mu )$. 
As for a Markov triple, see~\cite{bgl-book}, and since $P_t$ is symmetric and sub-Markovian, we can extend $P_t$ to a semigroup of contractions 
on $\mathbb{L}^p (\mu)$ for $p \in [ 1, \infty ]$. 
 We also let $\mathcal{E} (f) : = \mathcal{E} (f,f)$ and
$$
\| f \|_{\mathcal{E}}^2 : = \| f \|_{\mathbb{L}^2 (\mu)}^2 + \mathcal{E} (f) 
$$
for $f \in \mathcal{D} (\mathcal{E})$. 
We assume that $( \X , \tau , \mu, \mathcal{E} )$ is a Riemannian 
energy measure space in the sense of~\cite[Def.~3.16]{AGS_BE}, denoted $REM$ in this work. 
A basic example of a REM space is a Riemannian Markov triple as in Section~\ref{sec-riemannian}.
In this case, $( \mathcal{E} , \mathcal{D} (\mathcal{E}) )$ is canonically defined 
by completion of $( f, f ) \mapsto \int | \nabla f |^2 \, d \mu$. 
\textsf{RCD} spaces introduced in \cite{AGMR,AGS14} 
are another important class of $REM$ spaces. 
In this case, $\mathcal{E}/2$ is given by the $\mathbb{L}^2$-Cheeger energy functional. 

To make this presentation concise, we prefer to state the crucial properties of a $REM$ space
instead of its precise definition. Indeed the definition consists in several notions, which will be used 
only indirectly through these properties: 

\begin{itemize}
\item
The intrinsic distance $d_ {\mathcal{E}}$ associated with $( \mathcal{E} , \mathcal{D} ( \mathcal{E} ) )$, in the sense of \cite[Sec.~3.3]{AGS_BE},
becomes a distance function, further denoted $d$. 
It is compatible with the topology $\tau$ 
and the space $(\X, d)$ is complete \cite[Def.~3.6]{AGS_BE} and length metric \cite[Thm.~3.10]{AGS_BE}. 
\end{itemize}
We let $\mathrm{Lip}_b (\X)$ denote the set of bounded Lipschitz functions on $\X$ 
(with respect to $d$). 
Let $|\nabla f| : \X\to \mathbb{R}$ be 
the local Lipschitz constant of a Lipschitz function $f$ on $\X$: 
\[
| \nabla f | (x) 
: = 
\limsup_{y \to x} \frac{ | f (y) - f (x) | }{ d (x,y) } \cdot 
\]
\begin{itemize}
\item
$\mathcal{E}/2$ coincides with 
the $\mathbb{L}^2$-Cheeger energy associated with $d$, defined for $f\in \mathbb{L}^2 (\mu)$ by
$${\rm Ch}(f):=\inf\left\{\liminf_{n\to\infty}\frac 12\int |\nabla f_n|^2d\mu\,;\, f_n\in \mathrm{Lip}_b (\X),\, f_n\to f\, \mathrm{in} \, \mathbb{L}^2(\mu)\right\}.$$
As a result, $( \mathcal{E} , \mathcal{D} ( \mathcal{E} ) )$ 
admits a carr\'e du champ, i.e.\ there is a symmetric bilinear map 
$\Gamma : 
\mathcal{D} ( \mathcal{E} ) \times \mathcal{D} ( \mathcal{E} ) 
\to \mathbb{L}^1 ( \mu )
$ 
such that 
\[
\mathcal{E} ( f, g ) = \int \Gamma ( f , g ) \, d \mu. 
\]
As on smooth spaces, $L$ and $\Gamma$ satisfy the diffusion property 
\eqref{eq:diffusion2}. 
The coincidence of $\mathcal{E}/2$ and 
the Cheeger energy makes many connections between $d$ and $\Gamma$. 
For instance, 
$\mathcal{D} ( \mathcal{E} ) \cap \mathrm{Lip}_b (\X)$ 
is dense in $\mathcal{D} (\mathcal{E})$ with respect to 
$\| \cdot \|_{\mathcal{E}}$. 
In addition,
\begin{equation} \label{eq:Gamma-nabla}
\Gamma (f) \leq | \nabla f |^2 \quad \mbox{$\mu$-a.e.}
\end{equation} 
for any Lipschitz $f \in \mathcal{D} (\mathcal{E})$. 
See \cite[Thm.~3.12]{AGS_BE} and \cite[Thm.~3.14]{AGS_BE} for all these facts. 
\end{itemize}
Note that $\mathcal{D} ( \mathcal{E} ) \cap \mathbb{L}^\infty ( \mu )$ 
is an algebra and $\Gamma$ satisfies the Leibniz rule: 
\[
\Gamma ( f g , h ) = f \Gamma ( g , h ) + g \Gamma ( f , h )
\quad 
\mbox{for 
$f,g \in \mathcal{D} ( \mathcal{E} ) \cap \mathbb{L}^\infty ( \mu )$
and 
$h \in \mathcal{D} (\mathcal{E})$.
}
\] 

We state further assumptions for our main theorem. 
Fix a reference point $o \in \X$. 

\begin{eass} 
~
\begin{itemize}
\item[{\rm(Reg1)}]
There is $\alpha_0 > 0$ such that~\eqref{eq:integrable} holds. 
\item[{\rm(Reg2)}]
$( \X, \tau )$ is locally compact. 
\end{itemize}
\end{eass}

Assumption~(Reg1) 
is equivalent to the condition (MD.exp) in \cite{AGS_BE} 
(see e.g.\  the comments after Equation (3.13) in \cite{AGS_BE}). 
This integrability condition yields the conservativity of $P_t$,~i.e.
\[
\int P_t f \, d \mu = \int f \, d \mu 
\]
for $f \in \mathbb{L}^1 (\mu)$ (see~\cite[Thm.~3.14]{AGS_BE}). 
This is equivalent to $P_t 1 = 1$ $\mu$-a.e, that is, the semigroup is Markovian (instead of sub-Markovian).
In fact~\eqref{eq:integrable} is a nearly optimal condition 
to ensure that the semigroup is conservative 
(see \cite[Rmk.~4.21]{AGS13}). Thus it is not restrictive. 

Assumption~(Reg2) implies 
that any closed bounded set in $\X$ is compact 
(see e.g.~\cite[Prop.~2.5.22]{BBI}). 
Moreover, $(\X,d)$ is a geodesic space 
(see e.g.~\cite[Thm.~2.5.23]{BBI}). 
As a result, $(\mathcal{P}_2 (\X), W_2 )$ is also a geodesic space 
(see e.g.~\cite[Cor.~1 and Prop.~1]{Lisini:2006be}). 

\smallskip

In this framework, we should be careful when defining the operator  $\Gamma_2$ in~\eqref{eq:Gamma2} 
since $\Gamma (f)$ may not belong to $\mathcal{D} ( L )$ 
even for a sufficiently nice $f$. 
To avoid such a technical difficulty, and following~\cite[Def.~2.4]{AGS_BE}, we employ a weak form of the $CD (R,m)$ condition :
\begin{edefi}[Weak $CD(R,m)$ condition]
\label{def-weak-cd}
Let $R\in\R$ and $m>0$.  
 We say that the $REM$ space $(\X , \tau , \mu , \mathcal{E} )$ satisfies  
a weak  $CD (R,m)$ condition if, for all $f \in \mathcal{D} (L)$ 
with $L f \in \mathcal{D} (\mathcal{E})$ and all 
$g \in \mathcal{D} (L) \cap \mathbb{L}^\infty (\mu)$ 
with $g \geq 0$ and $L g \in \mathbb{L}^\infty (\mu)$, 
\begin{equation} \label{eq:weak-BE}
\frac12 \int \Gamma (f) L g \, d \mu 
- \int \Gamma ( f , L f ) g \, d \mu 
\geq R \int \Gamma (f) g \, d \mu 
+ \frac{1}{m} \int ( L f )^2 g \, d \mu . 
\end{equation}
\end{edefi}
Now we are ready to state our main theorem in this framework. 
%%%%
\begin{ethm} 
\label{thm:main_mms}
Let $( \X , \tau , \mu , \mathcal{E} )$ be a Riemannian energy measure space 
satisfying the above regularity assumptions (Reg1) and (Reg2). Let   $R \in \mathbb{R}$ and 
 $m > 0$. 

 If inequality~\eqref{eq-contraction-square-general-2} holds for any $t \geq 0$ and probability densities $f,g \in \mathbb{L}^1 (\mu)$ 
with $f \mu , g \mu \in \mathcal{P}_2 (\X)$, then the weak $CD ( R, m )$ condition of Definition~\ref{def-weak-cd} holds. 
In particular, the conditions (ii) and (iii) in Theorem~\ref{thm-legros} 
are equivalent to the weak $CD (R,m)$ condition.
\end{ethm}

Note that~\eqref{eq-contraction-square-general-2} 
yields a $W_2$-contraction 
\begin{equation}
\label{eq-contraction-easy2}
W_2^2(P_tf d\mu,P_tg d\mu)\leq e^{-2Rt} W_2^2(f d\mu,g d\mu)
\end{equation}
by neglecting the term involving $m$.
Then, by \cite[Cor.~3.18]{AGS_BE}, \eqref{eq-contraction-easy2} implies a $CD (R,\infty)$ condition 
in the sense of~\eqref{eq:weak-BE}. 
This fact is very helpful for further discussion in the sequel 
since it ensures regularity of the space in many respects. 
As a regularization property of $P_t$, we have 
\begin{equation} \label{eq:Lip-infinity}
\mbox{
$P_t f \in \mathrm{Lip}_b (\X)$  
for $f \in \mathbb{L}^2 (\mu) \cap \mathbb{L}^\infty (\mu), \; t>0$ 
}
\end{equation}
(see \cite[Thm.~3.17]{AGS_BE}; More precisely, 
$P_t f$ has a version which belongs to $\mathrm{Lip}_b (\X)$). 
In addition, $(\X, d, \mu )$ becomes an $\mathsf{RCD} (R, \infty)$ space 
(see \cite[Thm.~4.17]{AGS_BE}). 
Then, for a probability density $f$ with respect to $\mu$, 
$( ( P_t f ) \mu )_{t \geq 0}$ is a gradient flow of $\mathrm{Ent}_\mu$ 
in the sense of the $R$-evolution variational inequality \cite[Thm.~6.1]{AGMR}. 
As a consequence, we obtain the following properties: 
\begin{itemize}
\item\label{page-dual}
We can extend the action of $P_t$ to $\nu \in \mathcal{P}_2 (\X)$
in the sense that $P_t \nu$ is a solution to the $R$-evolution 
variational inequality and that $P_t \nu = ( P_t f ) \mu$ if $\nu = f \mu$. 
In particular, $( P_t \nu )_{t \geq 0}$ becomes a continuous curve 
in $( \mathcal{P}_2 (\X), W_2 )$, see \cite[Thm.~6.1]{AGMR}. 
In addition, $\nu \mapsto P_t \nu$ is a continuous map from 
$( \mathcal{P}_2 (\X), W_2 )$ to itself, see \cite[Eq. (7.2)]{AGMR}.

\item
$P_t \nu \ll \mu$ for $\nu \in \mathcal{P}_2 (\X)$ and $t > 0$, and its density $\rho_t$ satisfies 
$\ent{\mu}{\rho_t} \in \mathbb{R}$. 
This property is included 
in the definition of the $R$-evolution variational inequality, see e.g.\ \cite[Def.~2.5]{AGMR}. 
Recall that,  under~\eqref{eq:integrable}, $\ent{\mu}{\rho}$ is well-defined 
and $\ent{\mu}{\rho} \in ( - \infty , \infty ]$ 
for $\rho : \X \to [ 0, \infty ]$ 
with $\rho \mu \in \mathcal{P}_2 (\X)$, see e.g.~\cite[Sec.~7]{AGS13}. 

\item
There is a positive symmetric measurable function $p_t (x,y)$ such that 
$P_t$ coincides with the integral operator associated with $p_t$, see \cite[Thm.~7.1]{AGMR}. 

\item 
For any bounded measurable $h$ and $\nu \in \mathcal{P}_2 (\X)$, 
we have 
\begin{equation}
\label{eq-dual}
\int h \, d P_t \nu = \int P_t h \, d \nu, 
\end{equation}
see \cite[Prop.~3.2]{AGS_BE}. 
By the monotone convergence theorem, we can extend this identity 
to those $h$ which are bounded only from below (or above). 
\item For any $f \in \mathcal{D} (L)$ and $h \in \mathcal D({\mathcal E})$
we have the integration by parts formula
\begin{equation}
\label{eq-ipp-mms}
 \int \Gamma (h , f) \, d \mu =- \int h \, Lf \, d\mu.
\end{equation}
\end{itemize}

%%%%%%%%%%
\subsection{Proof of (ii) $\Rightarrow$ (iii) in Theorem~\ref{thm-legros}} \label{subsec:3-4}

As announced, before entering the proof of Theorem~\ref{thm:main_mms}, we first complete the proof of (ii) $\Rightarrow$ (iii) in
Theorem~\ref{thm-legros} with the aid of preparations in Section~\ref{subsec:frame_mms}. 

We first check that~\eqref{eq-contraction-sh} yields~\eqref{eq-contraction-easy2}. 
As we did in Section~\ref{subsec:frame_mms}, by using~\eqref{eq-contraction-square-general-2} 
instead of~\eqref{eq-contraction-sh},~we~get
\begin{equation} \label{eq:s-contraction}
s_{\frac Rm} \left( 
    \frac12 W_2 ( P_t f \mu , P_t g \mu )
\right)^2
\leq 
e^{-2Rt}\,
s_{\frac {R}{m}} \left(
    \frac {1}{2} W_2 ( f \mu , g \mu )
\right)^2 
\end{equation}
by neglecting the term involving $m$.
From this inequality, we can extend $P_t$ to a map 
from $\mathcal{P}_2 (\X)$ to itself, in a canonical way. 
Moreover, in~\eqref{eq:s-contraction}
we can replace $f \mu$ and $g \mu$ 
with any $\nu_0 , \nu_1 \in \mathcal{P}_2 (\X)$ respectively. 
Then we obtain~\eqref{eq-contraction-easy2} by a similar argument as in Proposition~\ref{prop-3-4}.
Thus, as discussed in Section~\ref{subsec:frame_mms},
$(\X , d, \mu )$ is an $\mathsf{RCD} ( R, \infty )$ space and 
all properties at the end of Section~\ref{subsec:frame_mms} become available. 
We remark that the extension of $P_t$ given on the basis of 
\eqref{eq:s-contraction} coincides with the one given 
by the $\mathsf{RCD} ( R, \infty )$ property. 

In Section~\ref{sec-proofdebut}, we already pointed out that 
we only need to show that $P_t$ fulfills all the assumptions for $\varphi_t$ in 
Proposition~\ref{prop-3-4} with $(Y , d_Y ) = ( \mathcal{P}_2 (\X) , W_2 )$ and 
$U = \entf{\mu}$. Here we are extending the definition of $\entf{\mu}$ 
so that, for $\nu \in \mathcal{P}_2 (\X)$, 
$\ent{\mu}{\nu} = \ent{\mu}{{ d \nu }/{ d \mu } }$ if $\nu \ll \mu$ 
and $\ent{\mu}{\nu} = \infty$ otherwise. 
By taking observations at the beginning of this section into account, 
it suffices to prove that~\eqref{eq-contraction-sh} implies 
\begin{multline*}
s_{\frac Rm} \left( 
    \frac12 W_2 ( P_t \nu_0 , P_t \nu_1 )
\right)^2
\leq 
e^{-2Rt}\,
s_{\frac {R}{m}} \left(
    \frac {1}{2} W_2 ( \nu_0 , \nu_1 )
\right)^2 
\\
- \frac{1}{2 m} \int_0^t 
e^{-2R(t-u)} 
\PAR{ 
    \ent{\mu}{ P_u \nu_0 } - \ent{\mu}{ P_u \nu_1 }
} ^2
du
\end{multline*} 
for $\nu_0 , \nu_1 \in \mathcal{P}_2 (\X)$ and $t > 0$. 
But this is true since $P_\delta \nu_0 , P_\delta \nu_1 \ll \mu$ for any $\delta \in (0,t)$, 
so that 
\begin{multline*}
s_{\frac Rm} \left( 
    \frac12 W_2 ( P_t \nu_0 , P_t \nu_1 )
\right)^2
\leq 
e^{-2R ( t - \delta )}\,
s_{\frac {R}{m}} \left(
    \frac {1}{2} W_2 ( P_\delta \nu_0 , P_\delta \nu_1 )
\right)^2 
\\
- \frac{1}{2 m} \int_{\delta}^{t} 
e^{-2R(t-u)} 
\PAR{ 
    \ent{\mu}{ P_{u} \nu_0 } - \ent{\mu}{ P_{u} \nu_1 }
} ^2
du
\end{multline*} 
by~\eqref{eq-contraction-sh} and the bound $\sinh^2(x) \geq x^2$; moreover $P_\delta \nu_i \to \nu_i$ in $W_2$ as $\delta \downarrow 0$ for $i=0,1$: this gives
the assertion. 
Hence the proof of (ii) $\Rightarrow$ (iii) in Theorem~\ref{thm-legros} is completed
and thus it is sufficient to show the main assertion of Theorem~\ref{thm:main_mms} 
to complete the proof of our equivalence result. 

%%%%%%%%%
\subsection{Construction of the path $(\tilde{g}_s)_{s\geq 0}$}
\label{subsec-path}

In this section, we build the path $\tilde{g}_s$ mentioned in Section~\ref{sec-how},
 under~\eqref{eq-contraction-square-general-2}.
Recall that $(\X , d, \mu )$ is now an $\mathsf{RCD} (R, \infty)$ space 
as remarked at the end of Section~\ref{subsec:frame_mms}.
For $x \in \X$ and $r > 0$, we denote the open ball of 
radius $r$ centered at $x$ by $B_r (x)$.

For this we first define $g (=\tilde{g}_0)$. We take $g$ in a more tractable (but large enough) class than the full class of  Definition~\ref{def-weak-cd}.
Fix $\alpha > \alpha_0$ with $\alpha_0$ as in~\eqref{eq:integrable}, $\lambda \in ( 0 , 1 )$ and $g_0 : \X \to \mathbb{R}$ Lipschitz 
with compact support. 
Let us define $g$ as follows: 
\begin{equation} \label{eq:g} 
g := \frac{1}{Z} \left( 
    ( 1 - \lambda ) g_0 
    + 
    \lambda \exp ( - \alpha d ( x , o )^2 )
\right)
\end{equation}
where $Z > 0$ is a normalizing constant 
such that $g \mu \in \mathcal{P} (\X)$. 
Note that~\eqref{eq:integrable} yields $g \mu \in \mathcal{P}_2 (\X)$. 
We fix $g$ until the end of the proof of Proposition~\ref{prop:preCD} below. 
We can define the $\mathbb{L}^2$-Cheeger energy functional $\mathcal{E}_g/2$
associated with $d$ and the probability measure $g \mu$. Let $\mathcal{D} ( \mathcal{E}_g )$ be 
the set of $f \in \mathbb{L}^2 (g \mu)$ with $\mathcal{E}_g (f) < \infty$. 
Recall that $\mathcal{D} ( \mathcal{E}_g )$ is complete with respect 
to $\| \cdot \|_{\mathcal{E}_g}$. 

To define the path $(\tilde{g}_s)_{s\geq 0}$ we need the corresponding generator $L^g$, and for this
we show the following auxiliary lemma. 
\begin{elem} \label{lem:g-bilinear}
 In the above notation,
$\mathcal{D} ( \mathcal{E} ) \subset \mathcal{D} ( \mathcal{E}_g )$ 
and 
\begin{equation} \label{eq:E_g}
\mathcal{E}_g ( f ) = \int \Gamma (f) g \, d \mu
\end{equation}
for $f \in \mathcal{D} (\mathcal{E})$. 
In addition, $( \mathcal{E}_g , \mathcal{D} ( \mathcal{E}_g ) )$ is bilinear. 
\end{elem}

We do not know whether~\eqref{eq:E_g} is valid 
 for any $f \in \mathcal{D} ( \mathcal{E}_g )$. 
Thus we have to be careful when we apply the integration by parts formula 
\eqref{eq-ipp} for $L^g$. 

\begin{Proof}
The former assertion follows from \cite[Lem.~4.11]{AGS13}. 
For the latter assertion, take $f, \tilde{f} \in \mathcal{D} ( \mathcal{E}_g )$. 
For each $n \in \mathbb{N}$, take also $\chi_n \in \mathrm{Lip}_b (\X)$ with $0 \leq \chi_n \leq 1$, 
$\chi|_{B_{n} (o)} \equiv 1$ and $\chi|_{B_{n+1} (o)^c} \equiv 0$. 
 
Since, for each $n \in \mathbb{N}$, $g$ is bounded away from $0$ on $B_n (o),$  
we have $f_n : = f \chi_n \in \mathcal{D} ( \mathcal{E} )$ 
by the locality of the Cheeger energy, see \cite[Prop.~4.8 (b)]{AGS13} 
and \cite[Lem.~4.11]{AGS13}. 
Moreover, 
$( f_n )_{n \in \mathbb{N}}$ forms a Cauchy sequence 
with respect to $\| \cdot \|_{\mathcal{E}_g}$ 
and hence $\| f_n - f \|_{\mathcal{E}_g} \to 0$. 
By the same argument, 
we have $\| \tilde{f}_n - \tilde{f} \|_{\mathcal{E}_g} \to 0$
for $\tilde{f}_n := \tilde{f} \chi_n$. 
By~\eqref{eq:E_g}, and recalling that $\Gamma$ is symmetric bilinear, we have
\[
\mathcal{E}_g ( f_n + \tilde{f}_n ) 
+ 
\mathcal{E}_g ( f_n - \tilde{f}_n ) 
= 
2 \left( 
    \mathcal{E}_g (f_n)
    + 
    \mathcal{E}_g (\tilde{f}_n)
\right). 
\]
Therefore the conclusion holds by letting $n \to \infty$. 
\end{Proof}
By Lemma~\ref{lem:g-bilinear}, 
$( \mathcal{E}_g , \mathcal{D} ( \mathcal{E}_g ) )$ is a closed bilinear form 
on $\mathbb{L}^2 ( g \mu )$. Hence there are an associated 
$\mathbb{L}^2$-semigroup $P^g_t$ of symmetric linear contraction 
and its generator $L^g$. 
By \cite[Prop.~4.8 (b)]{AGS13}, $\mathcal{E}_g$ is sub-Markovian. 
Thus $P^g_t$ satisfies the maximum principle, i.e. $P^g_t f \leq c$ if $f \leq c$ 
for $f \in \mathbb{L}^2 (g \mu)$ and $c \in \mathbb{R}$.
In addition, $\mathrm{Lip}_b (\X) \cap \mathcal{D} ( \mathcal{E}_g )$ 
is dense in $\mathcal{D} ( \mathcal{E}_g )$ 
with respect to $\| \cdot \|_{\mathcal{E}_g}$. 
Note that we can define $P^g_t$ and $L^g$ without bilinearity of $\mathcal{E}_g$ 
(see \cite[Sec.~4]{AGS13} and references therein). 
However, then they can be nonlinear and the integration by parts formula~\eqref{eq-ipp} may not hold. 
\begin{elem} \label{lem:g_dom}
In the above notation,
\begin{enumerate}
\item
$g \in \mathcal{D} ( \mathcal{E} ) \cap \mathbb{L}^\infty ( \mu )$  
and $\log g \in \mathcal{D} ( \mathcal{E}_g )$. 

\item
$\mathcal{D} (L) \subset  \mathcal{D} (L^g)$. 
\end{enumerate}
\end{elem}

\begin{Proof}
(i) The first claim follows 
from~\eqref{eq:Gamma-nabla} and~\eqref{eq:integrable}. 
For the second one, note that
\[
\mathcal{E}_g ( \log g ) \leq \int | \nabla \log g |^2 g \, d \mu.  
\]
It is the integrated form of~\eqref{eq:Gamma-nabla}
for $\mathcal{E}_g$ instead of $\mathcal{E}$. 
Then the claim follows 
from~\eqref{eq:integrable}. 

(ii) Let $f\in\mathcal D(L)$ and  $h \in \mathcal{D} (\mathcal{E}_g)$. 
Take $h_n \in \mathrm{Lip}_b (\X) \cap \mathcal{D} ( \mathcal{E}_g )$ 
for $n \in \mathbb{N}$ 
such that $\| h_n - h \|_{\mathcal{E}_g} \! \to \! 0$. 
By a truncation argument used in the proof of 
Lemma~\ref{lem:g-bilinear}, we may assume that each $h_n$ is supported 
on a bounded set, without loss of generality. 
Then $h_n \in \mathcal{D} ( \mathcal{E} ) \cap \mathbb{L}^\infty (\mu)$ 
and hence $h_n g \in \mathcal{D} ( \mathcal{E} )$. 
Thus the Leibniz rule, the assertion (i),~\eqref{eq-ipp-mms} and~\eqref{eq:Gamma-nabla} imply
\begin{align*}
\left| 
    \int \Gamma ( h_n , f ) g \, d \mu 
\right| 
& = 
\left| 
    \int \Gamma ( h_n g , f ) \, d \mu 
    - 
    \int h_n \Gamma ( g , f ) \, d \mu
\right|
\\
& \leq  
\left| 
    \int h_n ( L f ) g \, d \mu 
\right| 
+ 
\left| 
    \int h_n  \Gamma ( \log g , f ) g \, d \mu 
\right| 
\\
& \leq 
\| h_n \|_{L^2 (g\mu)} 
\left( 
    \| g \|_\infty \| Lf \|_{\mathbb{L}^2 (\mu)} 
    + 
    \left\| \frac{ | \nabla g |^2 }{g} \right\|_\infty 
    \mathcal{E} (f)^{1/2}
\right). 
\end{align*}
The definition of $g$ yields 
$\| | \nabla g |^2 / g \|_\infty \! < \! \infty$. 
Thus there is $C>0$ independent of $h$ and $n$ such~that 
\begin{align*}
\left| 
\mathcal{E}_g ( h_n , f )
\right| 
\leq
C \| h_n \|_{\mathbb{L}^2 (g\mu)}. 
\end{align*}
Here we used Lemma~\ref{lem:g-bilinear}.
By letting $n \! \to \! \infty$, we can replace $h_n$ with $h$ in this inequality. 
Hence $f \in \mathcal{D} ( L^g )$ since $h$ is arbitrary in $\mathcal{D} (\mathcal{E}_g)$. 
\end{Proof}

\bigskip

We can now define the path $(\tilde{g}_s)_{s \geq 0}$. 
Let $f \in \mathcal{D} (L) \cap \mathrm{Lip}_b (\X)$ with $\| f \|_{\infty} \leq 1/4$. 
We fix $f$ until the end of the following section, and observe that $f \in \mathbb{L}^2 ( g \mu )$. 
Then we let
\begin{equation} \label{eq:gs}
\tilde{g}_s : = g ( 1 + f - P^g_s f ). 
\end{equation}
By the $\mathbb{L}^{\infty}$-bound on $f$ and the maximum principle for $P^g_s$, 
we have 
\begin{equation} \label{eq:g-gs}
\frac12 g \leq \tilde{g}_s \leq 2 g. 
\end{equation}

In what follows, 
we may assume without loss of generality that $L^g f$ is not identically $0$. 
For, by~\eqref{eq-ipp-mms} and Lemma~\ref{lem:g_dom}, 
\begin{equation} \label{eq:g-log}
\int Lf \, g \, d \mu 
 = 
- \int \Gamma ( f , g ) \, d \mu 
 = 
- \int \Gamma ( f , \log g ) g \, d \mu 
 = 
\int L^g f  \, \log g \, g\, d \mu. 
\end{equation}
Thus, if $L^g f$ is identically $0$, then $\int  L f \, g \, d \mu = 0$;
hence~\eqref{eq:preCD1} below holds in this specific case 
(without the next section)
since the $CD (R , \infty )$ condition holds 
on our $\mathsf{RCD} (R, \infty)$ space. 

%%%%%%%%%%%%%%%%%%%
\subsection{Three key estimates}
\label{sec-3-estimates}

The proof of Theorem~\ref{thm:main_mms} is based on~\eqref{eq:preCD1} in Proposition~\ref{prop:preCD} below. In turn, this bound is based on the three key estimates in Lemmas~\ref{lem:first},~\ref{lem:second} and~\ref{lem:Ent_conv1}, which in the manifold case of Section~\ref{sec-proof-riemannian} correspond to~\eqref{eq-first-estimate},~\eqref{eq-second-estimate} 
and~\eqref{eq-third-estimate}. The proofs are a bit different since we use $\tilde{g}_s$ instead of $g_s$.

The Hopf-Lax semigroup $(Q_s)_{s \geq 0}$ given by~\eqref{eq:Hopf-Lax} 
will again play a crucial role. Required properties for $Q_s$ in this framework are given in~\cite[Sec.~3]{AGS13} or \cite[Sec.~3]{AGS_Sob} for instance.

\medskip

We begin with the {\it first estimate}, corresponding to~\eqref{eq-first-estimate}:
\begin{elem}[First estimate] \label{lem:first}
\begin{equation*}
\liminf_{s\rightarrow 0} \frac{W_2^2(P_t\tilde{g}_s\mu,P_tg\mu)}{2s^2}
\geq 
-\frac{1}{2}\int P_t( | \nabla f |^2 ) g\, d\mu + \int \Gamma(f,P_tf)g \, d\mu. 
\end{equation*}
\end{elem}
\begin{Proof}
It suffices to prove an lower bound on the right-hand side of~\eqref{eq:first-estimate0}.
By a rearrangement, 
\begin{equation}\label{eq:estimate1-1}
\int \frac{Q_s f P_t\tilde{g}_s- f P_tg}{s} d\mu 
= 
 \int \frac{Q_s f - f }{s} P_t ( \tilde{g}_s - g ) \, d \mu 
 +
 \int \frac{Q_s f - f }{s} P_t g  \, d \mu 
 + 
\int  f \frac{P_t ( \tilde{g}_s - g )}{s} \, d \mu . 
\end{equation}
Since $g \mu \in \mathcal{P} (\X)$, 
the Cauchy-Schwarz inequality yields $s^{-1} ( \tilde{g}_s - g ) \to - g \, L^g f$ 
in $\mathbb{L}^1 (\mu)$. 
Thus the last term in~\eqref{eq:estimate1-1} 
converges to $- \int f P_t ( g L^g f ) \, d \mu$. 
By  Lemma~\ref{lem:g-bilinear}, and as in Section~\ref{sec-proof-riemannian}, 
this quantity is equal to the second term on the right-hand side of the assertion. 
Moreover, by the general bound~\eqref{minoqsf}, the first term on the right-hand side of~\eqref{eq:estimate1-1} goes to 0. 
Finally, by~\eqref{minoqsf} and the Lebesgue dominated convergence theorem we conclude on the second term as in the Riemannian case of Section~\ref{sec-proof-riemannian}. 
More precisely, we have 
\begin{align*}
\liminf_{s \to 0} & \int \frac{Q_s f(x)  - f (x)}{s} P_t g (x) \, \mu (dx) 
\\
& \geq
- \frac{1}{2} \limsup_{s \to 0} 
\int \sup_{ y \in B ( x, \sqrt{ 4s \| f \|_\infty } ) \setminus \{ x \} } 
\left( 
\frac{ f (y) - f (x) }{ d(x,y) } 
\right)^2 P_t g (x) \mu (dx) 
= 
- \frac12 
\int | \nabla f |^2 P_t g \, d \mu.
\end{align*}
Thus the assertion holds. 
\end{Proof}

Next lemma deals with the {\it second estimate} and corresponds to~\eqref{eq-second-estimate}. 
\begin{elem}[Second estimate] \label{lem:second}
$$
\limsup_{s\rightarrow 0}\frac{W_2^2(\tilde{g}_s\mu,g\mu)}{2s^2}
\leq
\frac{1}{2 ( 1 -2 \| f \|_\infty )} 
\int \Gamma(f) g \, d\mu. 
$$
\end{elem}
\begin{Proof}
Again, by the dual form~\eqref{eq-23}, we need to bound $\int Q_s \psi \, \tilde{g}_sd\mu-\int \psi gd\mu$ uniformly from above on the bounded Lipschitz functions $\psi$. We can assume that $\psi$ is moreover supported on a bounded set.  
Then the function 
$(s_1 , s_2 ) \mapsto \int Q_{s_1}( \psi ) \tilde{g}_{s_2} \, d \mu$ 
satisfies the assumption of~\cite[Lem.~4.3.4]{ambrosio-gigli-savare} 
since we have~\eqref{eq:g-gs} and $\| Q_{s_1} \psi \|_\infty \leq \| \psi \|_\infty$. Thus, instead of~\eqref{eq-hl}, we obtain
\begin{align*} 
\frac{d}{ds} \int Q_s(\psi) \tilde{g}_s \, d \mu 
& \leq
\frac{d}{ds} \left. \int Q_s(\psi) \tilde{g}_{s_0} \, d \mu \right|_{s_0 = s} 
 + 
\frac{d}{ds} \left. \int Q_{s_0} (\psi) \tilde{g}_{s} \, d \mu \right|_{s_0 = s} 
\\
& = 
\int \SBRA{ 
  -\frac{1}{2}| \nabla Q_s\psi |^2 ( 1 + f - P^g_s f ) 
  - Q_s\psi \; L^g P_s^g f
} g \, d\mu 
\end{align*}
for a.e.\ $s>0$. 
Here the equality follows from \cite[Thm.~3.6]{AGS_Sob}, the properties
$\| Q_s \psi \|_{Lip} < \infty$,  
$\| Q_s \psi \|_{\infty} < \infty$ 
and the Lebesgue dominated convergence theorem.
Note that $Q_s \psi \in \mathcal{D} ( \mathcal{E}_g )$
since $Q_s \psi$ is Lipschitz with a bounded support. 
Thus, by virtue of Lemma~\ref{lem:g-bilinear} and~\eqref{eq:Gamma-nabla}, 
\begin{equation*}
- \int Q_s\psi \; ( L^g P_s^g f ) g \, d\mu 
 = 
\mathcal{E}_g ( Q_s\psi , P_s^g f )
 \leq 
\sqrt{ \mathcal{E}_g ( Q_s\psi ) \mathcal{E}_g ( P_s^g f ) }
 \leq 
\sqrt{ 
  \int | \nabla Q_s \psi |^2 g \, d\mu  
  \; 
  \mathcal{E}_g ( P_s^g f )}. 
\end{equation*}
By combining this estimate with the last one, 
we obtain 
\[
\frac{d}{ds} \int Q_s(\psi) \tilde{g}_s \, d \mu 
\leq
\frac{1}{ 2 ( 1 - 2 \|f\|_\infty ) } 
\mathcal{E}_g ( P_s^g f )
\leq 
\frac{1}{ 2 ( 1 - 2 \|f\|_\infty ) } 
\mathcal{E}_g ( f )
=
\frac{1}{ 2 ( 1 - 2 \|f\|_\infty ) } 
\int \Gamma( f ) g \, d \mu. 
\]
Here the second inequality follows from the spectral decomposition for quadratic forms 
and the equality follows from Lemma~\ref{lem:g-bilinear} again 
since $f \in \mathcal{D} (L) \subset \mathcal{D} ( \mathcal{E} )$. 
Thus the conclusion follows by integrating this estimate, as in the proof of~\eqref{eq-second-estimate}.
\end{Proof}

For the {\it third estimate}, we still require some preparation. 
We call $\mathcal{C}_2 (\X)$ the set of continuous functions  $\psi$ on $\X$ 
for which there exists $C > 0$ such that $| \psi (x) | \leq C ( 1 + d ( o , x )^2 )$. 
For $\psi \in \mathcal{C}_2 (\X)$ and $\nu \in \mathcal{P}_2 (\X)$, 
we have $\psi \in \mathbb{L}^1 (\nu)$. 
By assumption on $g$, $\psi \in \mathbb{L}^p (g \mu )$
for any  $\psi \in \mathcal{C}_2 (\X)$ and $p \in [ 1, \infty )$. 
The following lemma ensures integrability properties 
required in the proof of Lemma~\ref{lem:Ent_conv1} below.
\begin{elem} \label{lem:integrable}
~
\begin{enumerate}
\item
Let $J : = \{ \psi \in \mathbb{L}^2 (g \mu) \; | \; \psi g \mu \in \mathcal{P} (\X) \}$. 
Then $\psi g \mu \in \mathcal{P}_2 (\X)$ for any $\psi \in J$. 
Moreover, for $J_0 \subset J$ 
with $\sup_{\psi \in J_0} \| \psi \|_{\mathbb{L}^2 (g\mu)} < \infty$, 
we have 
\[
\sup_{\psi \in J_0} \int d (o,x)^2 \psi g \, d \mu < \infty. 
\]

\item
$\log P_u g \in \mathcal{C}_2 (\X)$ for $u \geq 0$. 
\end{enumerate}
\end{elem}
\begin{Proof}
(i) Using Assumption (Reg1) and~\eqref{eq:g}, this follows from
 \begin{align*}
\int d ( o , x )^2 \psi (x) g (x) \, \mu (dx)
\leq 
\left( \int d ( o , x )^4 g(x) \, \mu (dx) \right)^{1/2} 
\left( \int \psi^2 g \, d \mu \right)^{1/2}
< \infty.
\end{align*}

(ii) By~\eqref{eq:g} this is obvious for $u = 0$ and hence we consider the case $u > 0$.  
First of all, $\log P_u g$ is continuous on $\X$ since $P_u g > 0$.
Moreover, since $(\X, d, \mu )$ is an $\mathsf{RCD} (R,\infty)$ space, 
we have the log-Harnack inequality 
\begin{equation*} 
P_u ( \log g ) (o) - \frac{R d (x,o)^2}{2( e^{2Ru} - 1 )}
\leq \log P_u g (x) \leq \log \| g \|_{\infty} 
\end{equation*}
(see \cite[Lem.~4.6]{AGS_BE} or \cite[Prop.~4.1]{LiHQ15}). 
Moreover $\log g \in \mathcal{C}_2 (\X)$ and $P_u \delta_o \in \mathcal{P}_2 (\X)$ by the properties after~\eqref{eq:Lip-infinity}, so
we have $\int \log g \, d P_u \delta_o = P_u ( \log g ) (o) \in \R$. 
Thus $\log P_u g \in \mathcal{C}_2 (\X)$.
\end{Proof}

We recall characterizations of convergence in $W_2$ for later use. 
Let $\nu_n \in \mathcal{P}_2 (\X)$, $n \in \mathbb{N}$ and $\nu \in \mathcal{P}_2 (\X)$. 
Then $W_2 ( \nu_n , \nu ) \to 0$ is equivalent to either of the following 
(see e.g. \cite[Thm.~6.9]{villani-book2}): 
\begin{itemize}
\item
$\nu_n \to \nu$ weakly and 
$\displaystyle 
\sup_{n \in \N} \int d ( o, x )^2 \nu_n (dx) < \infty
$, 

\item
$\displaystyle \lim_{n \to \infty} \int \psi \, d \nu_n = \int \psi \, d \nu$
for any $\psi \in \mathcal{C}_2 (\X)$. 
\end{itemize}

We now turn to the {\it third estimate}. 
\begin{elem}[Third estimate] \label{lem:Ent_conv1}
$$
\liminf_{s \to 0} 
\frac{1}{s^2}
\int_0^t \! e^{-2R(t-u)}\SBRA{\ent{\mu}{P_u \tilde{g}_s}-\ent{\mu}{P_u g}}^2du
\geq
\int_0^t \! e^{-2R(t-u)} \SBRA{ \int P_u \big( g L^g f \big) \log P_u g \, d \mu }^2 du . 
$$
\end{elem}

\begin{Proof}
By the Fatou lemma, it suffices to show 
\[
\liminf_{s \to 0} \SBRA{ \frac{ \ent{\mu}{P_u \tilde{g}_s} - \ent{\mu}{P_u g} }{s} }^2 
\geq 
\SBRA{ \int P_u \big( g L^g f \big) \log P_u g \, d \mu }^2  
\]
for each $u > 0$. 
By~\eqref{eq:g-gs} 
and since
$\ent{\mu}{P_u g} \in \R$, 
we have $P_u \tilde{g}_s \log P_u g, P_u g \log P_u g \in \mathbb{L}^1 (\mu)$. Moreover
$a^2 \geq (a+b)^2/(1+\delta)-b^2 / \delta$ for $\delta > 0$ and 
$$
0\leq x\log x-x+1\leq (x-1)^2
$$
for $x\geq 0$, 
so
$$
\left( \ent{\mu}{P_u \tilde{g}_s} - \ent{\mu}{P_u g} \right)^2 
\geq 
\frac{1}{1 + \delta } 
\left( 
    \int (P_u \tilde{g}_s - P_u g ) \log P_u g \, d \mu 
\right)^2 
- \frac{1}{\delta} 
\left( \int \frac{( P_u \tilde{g}_s - P_u g )^2}{P_u g} d\mu
\right)^2.
$$
By the Cauchy-Schwarz inequality for $P_u$, 
$$
\limsup_{s \to 0} \frac{1}{s} \!
\int \! \frac{ ( P_u \tilde{g}_s - P_u g )^2}{P_u g} d\mu
 \leq 
\limsup_{s \to 0}
\frac{1}{s} \!
\int \! P_u \! \left( 
    \frac{( \tilde{g}_s - g )^2}{g}
\right) 
d \mu 
 = 
\limsup_{s \to 0}\,\, s\!\! 
\int \! \left| \frac{ P_s^g f - f }{s} \right|^2 \! g \, d \mu 
= 0. 
$$
Since $\delta > 0$ is arbitrary, it suffices to show
\begin{equation} \label{eq:Ent_conv1}
\lim_{s \to 0} \frac{1}{s} \int P_u \big( g( P^g_s f - f ) \big) \log P_u g \, d \mu 
= 
\int P_u \big( g L^g f \big) \log P_u g \, d \mu  
\end{equation}
in order to complete the proof. 
Here the well-definedness of the right-hand side 
is included in the assertion. 
Since $r \mapsto r_+$ is 1-Lipschitz, 
$s^{-1} ( P^g_s f - f )_{+} = ( s^{-1} ( P^g_s f - f ) )_{+}$ 
converges to $(L^g f )_{+}$ 
in $\mathbb{L}^2 (g \mu)$ and hence in $\mathbb{L}^1 ( g \mu )$. 
By~\cite[Thm.~4.16 (d)]{AGS13}, $\int L^g f \, g \, d \mu = 0$. 
Hence $\|( L^g f )_{+} \|_{\mathbb{L}^1 (g\mu)} >0$ since $L^g f$ is not identically $0$ 
(as assumed at the end of Section~\ref{subsec-path}). 
Thus $\| (P^g_s f - f)_{+} \|_{\mathbb{L}^1 ( g \mu )} > 0$
for sufficiently small $s > 0$. 
Let us now define $\nu^f_s , \nu^f_0 \in \mathcal{P} (\X)$ as follows: 
\begin{align*}
\nu^f_s : = 
\frac{ ( P^g _s f - f )_{+} }{ \| ( P^g_s f - f )_{+} \|_{\mathbb{L}^1 (g \mu)} } g \mu,  
\qquad 
\nu^f_0 : = 
\frac{( L^g f )_+}{\| ( L^g f )_+ \|_{\mathbb{L}^1 (g \mu)} } g \mu . 
\end{align*} 
Then $\nu^f_s \to \nu_0^f$ weakly in $\mathcal{P} (\X)$ as $s \to 0$. 
Moreover, $\nu^f_s \in \mathcal{P}_2 (\X)$ for $s \geq 0$ by (i) in Lemma~\ref{lem:integrable}
since $f, P_s^g f, L^g f \in \mathbb{L}^2 (g \mu)$, and 
$W_2 ( \nu^f_s , \nu^f_0 ) \to 0$ as $s \to 0$, 
again by (i) in Lemma~\ref{lem:integrable} applied with 
$J_0 =\{\| ( P^g_s f - f )_+ \|_{\mathbb{L}^1 (g \mu )}^{-1}(P_s^g f - f)_+, s>0\}$, 
and the remark after it. 
Then, likewise, $P_u \nu^f_s \in \mathcal{P}_2 (\X)$ for $u,s \geq 0$ and 
\begin{equation} \label{eq:W2conv}
\lim_{s \to 0} W_2 ( P_u \nu^f_s , P_u \nu^f_0 ) = 0 
\end{equation}
by~\eqref{eq-contraction-easy2}. 
By Lemma~\ref{lem:integrable} again, $\log P_u g \in \mathcal{C}_2 (\X)$ 
and in particular $\log P_u g \in \mathbb{L}^1 ( P_u \nu^f_0 )$. 
Hence, by~\eqref{eq:W2conv} and the remark after Lemma~\ref{lem:integrable}, 
we obtain 
\begin{multline*}
\lim_{s \to 0} \frac{1}{s}  \int P_u ( g ( P_s^g f - f )_+ ) \log P_u g \, d \mu
= 
\lim_{s\to 0} \frac{ \| ( P^g_s f - f)_+ \|_{\mathbb{L}^1 (g\mu)} }{ s }
\int \log P_u g \, d P_u \nu^f_s 
\\
= 
\| ( L^g f )_+ \|_{\mathbb{L}^1 (g\mu)}
\int \log P_u g \, d P_u \nu^f_0 
= 
\int P_u ( g ( L^g f )_+ ) \log P_u g \, d \mu  \, \in \R.
\end{multline*}
We can apply the same argument to $( P^g_s f - f )_-$ instead of 
$( P^g_s f - f )_+$ to show the corresponding assertion. 
In particular, the integral in the right-hand side of~\eqref{eq:Ent_conv1} 
is well-defined and these two claims yield~\eqref{eq:Ent_conv1}. 
\end{Proof}

%%%%%%%%%%%%%%%%%%%%%%%%%%
\subsection{Conclusion of the proof of Theorem~\ref{thm:main_mms}}
\label{subsec:pf_mms}
Let $g$ be as in the last section, that is, given by~\eqref{eq:g}.  To proceed, we recall the notion of semigroup mollification 
introduced in \cite[Sec.~2.1]{AGS_BE}. 
Let $\kappa \in \mathcal{C}^\infty_c (( 0, \infty ))$ with $\kappa \geq 0$ 
and $\int_0^\infty \kappa (r) \, d r = 1$. 
For $\ep > 0$ and $f \in \mathbb{L}^p (\mu)$ with $p \in [ 1 , \infty ]$, 
we define $\mathfrak{h}_\ep f$ by 
\[
\mathfrak{h}_{\ep} f 
: = 
\frac{1}{\ep}  \int_0^\infty 
P_r f \; \kappa \left( \frac{r}{\ep} \right) 
\, d r. 
\]
It is immediate that 
$\| \mathfrak{h}_\ep f - f \|_{\mathcal{E}} \to 0$ 
as $\ep \to 0$ for 
$f \in \mathcal{D} ( \mathcal{E} )$. 
Moreover, for $f \in \mathbb{L}^2 ( \mu ) \cap \mathbb{L}^\infty ( \mu )$, 
$\mathfrak{h}_\ep f , L (\mathfrak{h}_\ep f) \in \mathcal{D} (L) \cap \mathrm{Lip}_b (\X)$. 
Here the latter one comes from the following representation: 
\[
L \mathfrak{h}_\ep f 
= 
- \frac{1}{\ep^2} \int_0^\infty 
  P_r f \; \kappa' \left( \frac{r}{\ep} \right) 
\, d r .
\] 

\begin{eprop} \label{prop:preCD} Following the same assumptions as in Theorem~\ref{thm:main_mms}, let $f = \mathfrak{h}_\e f_0$ for some $\e > 0$ 
and $f_0 \in \mathbb{L}^2 (\mu) \cap \mathbb{L}^\infty (\mu)$. 
Then $\Gamma (f) \in \mathcal{D} ( \mathcal{E} )$, and for $g$ as above 
\begin{equation} \label{eq:preCD1}
\frac12 \int \Gamma ( \Gamma (f) , g ) \, d\mu 
+
\int \Gamma ( f , L f ) g \, d \mu
\leq 
- R \int \Gamma(f) g \, d\mu 
-
\frac{1}{m} \PAR{ \int Lf \, g \, d\mu }^2 . 
\end{equation}
\end{eprop}
\begin{Proof}
By assumption, $f \in \mathcal{D} (L) \cap \mathrm{Lip}_b (\X)$.
Moreover, $\Gamma (f) = | \nabla f |^2$ 
$\mu$-a.e.\ by \cite[Thm.~3.17]{AGS_BE}. 
Let $\eta > 0$ be so small that $\eta \| f \|_\infty \leq 1/4$. 
By applying Lemma~\ref{lem:first}, Lemma~\ref{lem:second}
and Lemma~\ref{lem:Ent_conv1} to $\eta f$ instead of $f$ 
in~\eqref{eq-contraction-square-general-2}, 
\begin{multline*} 
-\frac{\eta^2}{2} 
\int  P_t \Gamma (f) \; g \, d\mu
+
\eta^2 \int \Gamma ( f ,P_t f ) g \, d\mu
\\
\leq
\frac{e^{-2Rt} \eta^2 }{2 ( 1 - 2 \eta \| f \|_\infty ) } 
\int \Gamma ( f ) g \, d\mu
- 
\frac{\eta^2}{m} \int_0^t e^{-2R(t-u)}
\PAR{ 
  \int
  P_u ( ( L^g f ) g ) \log P_u g 
  \, d\mu
}^2 
du. 
\end{multline*}
By dividing this inequality by $\eta^2$ and letting $\eta \to 0$, 
\begin{multline} \label{eq:preCD1-1}
-\frac{1}{2} 
\int  P_t \Gamma (f) \; g \, d\mu
+
\int \Gamma ( f ,P_t f ) g \, d\mu
\\
\leq
\frac{e^{-2Rt} }{2}
\int \Gamma ( f ) g \, d\mu
- 
\frac{1}{m} \int_0^t e^{-2R(t-u)}
\PAR{ 
  \int
  P_u ( ( L^g f ) g ) \log P_u g 
  \, d\mu
}^2 
du. 
\end{multline}
By virtue of mollification by $\mathfrak{h}_\e$, 
we have $L f \in \mathcal{D} (\mathcal{E})$ and 
\begin{align*}
\left. \frac{d}{dt} \right|_{t=0} 
\int \Gamma ( f ,P_t f ) g \, d\mu
& = 
- \frac{1}{\e^2} \int_0^\infty \kappa' \PAR{ \frac{r}{\e} }
\int \Gamma ( f , P_r f_0 ) g \, d\mu
d r 
 = 
\int \Gamma ( f , L f ) g \, d\mu .
\end{align*}
Note that $\Gamma (f) \in \mathcal{D} ( \mathcal{E} )$ 
 (hence the left-hand side of~\eqref{eq:preCD1} is well-defined).
This fact follows from \cite[Lem.~3.2]{Savare:2014jm} 
with the aid of mollification by $\mathfrak{h}_\ep$. 
Then, by Lemma~\ref{lem:Ent_conv2} below, 
we can differentiate~\eqref{eq:preCD1-1} at $t = 0$ to obtain 
\begin{align*} 
\frac{1}{2} 
\int \Gamma ( \Gamma (f), g ) \, d\mu
 +
\int \Gamma ( f , L f ) g \, d\mu
& \leq
- R \int \Gamma ( f ) g \, d\mu
- \frac{1}{m} 
\PAR{ 
  \int
  ( L^g f ) g \log g 
  \, d\mu
}^2 
\\ 
& = 
- R \int \Gamma ( f ) g \, d\mu
- \frac{1}{m} 
\PAR{ 
  \int
  ( L f ) g 
  \, d\mu
}^2 . 
\end{align*}
Here we have used~\eqref{eq:g-log} also in the last equality. 
This is nothing but the desired inequality. 
\end{Proof}
\begin{elem} \label{lem:Ent_conv2}
For $\psi \in L^2 ( g \mu )$, 
\[
\lim_{u \to 0} \int P_u ( \psi g ) \log P_u g \, d \mu 
= 
\int \psi g \log g \,d \mu . 
\]
\end{elem}
\begin{Proof}
We may assume $\psi \ge 0$ and $\psi g \mu \in \mathcal{P} (\X)$ 
without loss of generality. 
Then in particular $\psi g \mu \in \mathcal{P}_2 (\X)$ by Lemma~\ref{lem:integrable}.
First of all, 
\[
\int P_u ( \psi g ) | \log P_u g | \, d \mu < \infty 
\]
by a similar argument as in Lemma~\ref{lem:integrable}. Thus 
\[
\int P_{u} ( \psi g ) \log P_{u} g \, d \mu 
= \int \psi g P_{u} ( \log P_{u} g ) \, d \mu
\leq 
\int \psi g \log P_{2u} g \, d \mu 
\]
by the Fubini theorem and the Jensen inequality for $P_u$ 
 as integral operator. 
Now, for each $x$, $\lim_{u \to 0} W_2 ( P_u \delta_x , \delta_x ) = 0$ by the remark after Theorem~\ref{thm:main_mms}, and $g$ is bounded and continuous, so $P_u g(x) = \int g d P_u \delta_x \to g(x).$ Moreover $\log P_{2u} g\leq \log ||g||_\infty$ and $\psi g\mu$ is a probability measure, so
by the Fatou lemma
\begin{equation} \label{eq:Ent_conv2up}
\limsup_{u \to 0} 
\int P_{u} ( \psi g ) \log P_{u} g \, d \mu 
\leq 
\int \psi g \log g \, d \mu . 
\end{equation}

For the opposite bound, again by the Jensen inequality for $P_u$, 
\[
\int P_{u} ( \psi g ) \log P_{u} g \, d \mu 
\geq 
\int P_u ( \psi g ) P_u ( \log g ) \, d \mu 
= 
\int \log g P_{2u} ( \psi g ) \, d \mu . 
\]
Moreover $\log g$ is in $\mathcal{C}_2 (\X)$ and $W_2 ( P_{2 u} ( \psi g ) \mu , \psi g \mu ) \to 0$ as $u \to 0$, again by the remark after Theorem~\ref{thm:main_mms}. Hence, by the remark after Lemma~\ref{lem:integrable},
we obtain 
\begin{equation} \label{eq:Ent_conv2down}
\liminf_{u \to 0} \int P_{u} ( \psi g ) \log P_{u} g \, d \mu 
\geq 
\lim_{u \to 0} \int P_{2u} ( \psi g ) \log g \, d \mu 
= 
\int \psi g \log g \, d \mu . 
\end{equation}
Hence the conclusion follows from the combination of 
\eqref{eq:Ent_conv2up} and~\eqref{eq:Ent_conv2down}.
\end{Proof}

Now we are in turn to complete the proof of Theorem~\ref{thm:main_mms}. 
\begin{tProof}{Theorem~\ref{thm:main_mms}} The last crucial step consists in transforming 
$( 
  \int
  ( L f ) g 
  \, d\mu
)^2$
into
$
  \int
  ( L f )^2 g 
  \, d\mu
$
which will be done by a localization procedure. Let $f$ be as in Proposition~\ref{prop:preCD}.

Remark first that, by letting $\lambda \to 0$ in the definition~\eqref{eq:g}, we obtain 
\eqref{eq:preCD1} for $g_0$ instead of the function $g$ of~\eqref{eq:g}. 
To put the square inside the integral in~\eqref{eq:preCD1}, we need to \emph{localize} this inequality, and thus we employ a partition of unity. 
Let $\eta > 0$. Since $L f \in \mathrm{Lip}_b (\X)$, 
we can take $\delta > 0$ sufficiently small 
so that $| L f (x) - L f (y) | < \eta$ 
for any $x,y \in \X$ with $d(x,y) < 4 \delta$. 
Since $\mathrm{supp}\, g_0$ is compact, 
there is $\{ x_i \}_{i=1}^n \subset \X$ such that 
$\mathrm{supp}\, g_0 \subset \bigcup_{i=1}^n B_{\delta} (x_i)$
(note that we require the regularity assumption~(Reg2) 
only at this point). 
Let us define $\tilde{\psi}_i$ 
($i=1, \ldots , n$) 
by $\tilde{\psi}_i (x) : = 0 \vee ( 2 \delta - d (x_i , x ) )$ 
and 
\[
\psi_i (x) := 
\begin{cases}
\displaystyle 
\frac{ \tilde{\psi}_i (x) }{ \sum_{j=1}^n \tilde{\psi}_j (x) }
& 
\mathrm{if} \quad \tilde{\psi}_i (x) \neq 0 , 
\\
0 & \mathrm{if} \quad \tilde{\psi}_i (x) = 0. 
\end{cases}
\] 
Then $\psi_i \in \mathrm{Lip} (\mathrm{supp}\ g_0)$, $0 \leq \psi_i \leq 1$,  
$\mathrm{supp}\, \psi_i \subset B_{2 \delta} (x_i)$ 
and $\sum_{i=1}^n \psi_i (x) = 1$ for $x \in \mathrm{supp}\, g_0$. 
By applying~\eqref{eq:preCD1} for $\psi_i g_0 / \| \psi_i g_0 \|_{\mathbb{L}^1 (\mu)}$ 
instead of $g_0$, we have 
\begin{align*}
\frac12 & \int \Gamma ( \Gamma ( f ) , g_0  ) \, d \mu 
+ 
\int \Gamma (f , L f ) g_0 \, d \mu 
= 
\sum_{i =1}^n
\left( 
    \frac12 \int \Gamma ( \Gamma ( f ) , \psi_i g_0 ) \, d \mu 
    + 
    \int \Gamma (f , L f ) \psi_i g_0 \, d \mu 
\right)
\\
& \leq 
- R \int \Gamma (f) g_0 \, d \mu
- \frac{1}{m}
\sum_{i=1}^n \frac{1}{ \| \psi_i g_0 \|_{\mathbb{L}^1 (\mu)}}
\left( 
    \int ( Lf ) \psi_i g_0 \, d \mu
\right)^2 . 
\end{align*}
By the choice of $\delta$ and $\{ \psi_i \}_{i=1}^n$, with $\eta<1$, 
\begin{align*}
\sum_{i=1}^n \frac{1}{ \| \psi_i g_0 \|_{\mathbb{L}^1 (\mu)}}
\left( 
    \int ( Lf ) \psi_i g_0 \, d \mu
\right)^2 
& \geq 
(1-\eta) \sum_{i=1}^n \| \psi_i g_0 \|_{\mathbb{L}^1 (\mu)} Lf (x_i)^2 
- \eta 
\\
& \geq  
(1-\eta)\int  ( Lf )^2 g_0 \, d\mu-\eta-2\eta(1-\eta)\|Lf\|_{\infty}
\end{align*}
By letting $\eta \to 0$, 
\[
- \frac12 \int \Gamma ( \Gamma ( f ) , g_0  ) \, d \mu 
- 
\int \Gamma (f , L f ) g_0 \, d \mu 
\geq 
R \int \Gamma (f) g_0 \, d \mu
+ 
\frac{1}{m} \int ( Lf )^2 g_0 \, d \mu .
\] 

Let now $g \in \mathcal{D} (L) \cap \mathbb{L}^\infty (\mu)$ 
with $g \geq 0$ and $L g \in \mathbb{L}^\infty (\mu)$, 
as in Theorem~\ref{thm:main_mms}. 
By virtue of mollification by $\mathfrak{h}_\ep$, 
\eqref{eq:Gamma-nabla} and~\eqref{eq:Lip-infinity}, 
we have 
$
\Gamma (f) , \Gamma ( f , L f ), ( Lf )^2 
\in 
\mathbb{L}^1 (\mu) \cap \mathbb{L}^\infty (\mu) 
$. 
Thus we can replace $g_0$ in the last inequality 
with $g_1 \in \mathrm{Lip}_b (\X) \cap \mathcal{D} ( \mathcal{E} )$, 
by a standard truncation argument.
Then we can replace $g_1$ with $g$ 
since $\mathcal{D} (\mathcal{E}) \cap \mathrm{Lip}_b (\X)$ 
is dense in $\mathcal{D} (\mathcal{E})$ 
with respect to $\| \cdot \|_{\mathcal{E}}$.

Finally, we remove the mollification $\mathfrak{h}_\ep$. 
Let $f \in \mathcal{D} ( L )$ with $L f \in \mathcal{D} (\mathcal{E})$
and $f_n : = (-n) \vee f \wedge n$. 
Then we have,  from the integration by parts formula~\eqref{eq-ipp-mms},
\[
\frac12 \int \Gamma ( \mathfrak{h}_\ep f_n ) L g \, d \mu 
- 
\int \Gamma ( \mathfrak{h}_\ep f_n , L \mathfrak{h}_\ep f_n ) g \, d \mu 
\geq
R \int \Gamma ( \mathfrak{h}_\ep f_n ) g \, d \mu
+ 
\frac{1}{m} \int ( L \mathfrak{h}_\ep f_n )^2 g \, d \mu .
\]
By virtue of mollification by $\mathfrak{h}_\ep$, 
$\| \mathfrak{h}_\ep f_n - \mathfrak{h}_\ep f \|_{\mathcal{E}} \to 0$
and $\| L \mathfrak{h}_\ep f_n - L \mathfrak{h}_\ep f \|_{\mathcal{E}} \to 0$ 
as $n \to \infty$. 
Thus we obtain~\eqref{eq:weak-BE} 
by letting $n \to \infty$ and $\ep \to 0$ after it, 
with taking $L \mathfrak{h}_\ep f = \mathfrak{h}_\ep L f$ into account. 
\end{tProof}

%%%%%
\section{Links with functional inequalities}

%%%
{\bf{A new proof of the entropy-energy inequality}}

\medskip
We now
% In this section we 
consider the case where $R>0$ and $\mu$ is a probability measure. It is classical that the $CD(R,m)$ condition implies the entropy-energy inequality  
\begin{equation} \label{eq:En-En}
\ent{\mu}{f}\leq \frac{m}{2}\log\left(1+\frac{1}{mR}I(f)\right)
\end{equation}
for any function $f$ such that $\int fd\mu=1$. Here $I(f)=\int\Gamma(f)/fd\mu$ is the Fisher information of~$f$. This inequality is given in~\cite[Thm.~6.8.1]{bgl-book} for instance, and also in \cite[Cor. 3.28]{EKS13} via the $(R,m)$-convexity of $\entf{\mu}$. 

Inequality~\eqref{eq:En-En} improves upon the standard non dimensional logarithmic Sobolev inequality $\ent{\mu}{f}\leq I(f)/ 2R,$ a consequence of the $CD(R, \infty)$ condition. It leads 
 for example to a sharp bound on the instantaneous creation of the entropy of the heat semigroup in $\mathcal P_2(\X)$, namely 
$$
\ent{\mu}{P_tf}\leq \frac m2 \log\frac{1}{1-e^{-2Rt}}
$$
for all $f$ and $t>0$. For similar bounds, see also \cite[Prop. 2.17]{EKS13} for a gradient flow argument starting from the $(R,m)$-convexity of $\entf{\mu}$, and~\cite[Prop. 3.1]{BGG15} for Fokker-Planck equations on $\mathbb R^m$ with $R$-convex potentials.  

\smallskip

The two approaches of \cite{bgl-book} and \cite{EKS13} are rather involved, and we now give a formal (and below rigorous) and direct way of recovering~\eqref{eq:En-En} from the contraction inequality~\eqref{eq-contraction-square-general-2}  in Theorem~\ref{thm-legros} (which is equivalent to the $CD(R,m)$ condition). The key point is the (formal) identity 
\begin{equation} \label{eq:speed}
\limsup_{\delta \downarrow 0}\frac{W_2^2(P_{\delta  + t } f\mu, P_t  f\mu)}{\delta^2}
= I( P_t f)
\end{equation}
(see e.g.~\cite[Equation~(26)]{ov00}) and the classical identity $\frac{d}{du}\ent{\mu}{P_uf}=-I(P_uf)$. 
Indeed, from inequality~\eqref{eq-contraction-square-general-2} and the Fatou Lemma, for any $0\leq s<t$, 
\begin{align*}
I ( P_t f ) 
 = 
\limsup_{\delta \downarrow 0} \frac{ W_2^2(P_{t+\delta} f\mu,P_{t} f \mu) }{\delta^2} 
& \leq 
e^{-2R(t-s)} 
\limsup_{\delta \downarrow 0} \frac{ W_2^2(P_{s+ \delta} f \mu, P_s f \mu) }{\delta^2}
\\ &\quad - \! \frac{2}{m} \!  \int_s^t \! e^{-2R(t-u)} 
\liminf_{\delta \downarrow 0} 
\left( 
    \frac{ \ent{\mu}{P_{u+\delta} f} - \ent{\mu}{P_u f} \!}{\delta}
 \right)^2 \! du
\\
& = 
e^{-2R(t-s)} I ( P_s f )
- \frac{2}{m}\int_s^{t}e^{-2R(t-u)} I ( P_u f )^2 du. 
\end{align*}
%\begin{align*}
%I ( P_t f ) 
%& = 
%\limsup_{\delta \downarrow 0} \frac{ W_2^2(P_{t+\delta} f\mu,P_{t} f \mu) }{\delta^2} 
%\\
%& \leq 
%e^{-2R(t-s)} 
%\limsup_{\delta \downarrow 0} \frac{ W_2^2(P_{s+ \delta} f \mu, P_s f \mu) }{\delta^2}
%\\ &\quad\quad\quad\quad\quad\quad\quad\quad\quad - \frac{2}{m}\int_s^t e^{-2R(t-u)} 
%\liminf_{\delta \downarrow 0} 
%\left( 
%    \frac{ \ent{\mu}{P_{u+\delta} f} - \ent{\mu}{P_u f} }{\delta}
%\right)^2 du
%\\
%& = 
%e^{-2R(t-s)} I ( P_s f )
%- \frac{2}{m}\int_s^{t}e^{-2R(t-u)} I ( P_u f )^2 du. 
%\end{align*}
This yields the differential inequality
$$
\frac{d}{dt}I(P_tf)\leq - 2R \, I(P_tf)-\frac2mI(P_tf)^2
$$
and then 
\begin{equation}\label{eq:I2}
I ( P_t f ) \le \frac{ m R I (f) }{ e^{2Rt} ( I(f) + m R ) - I (f) } 
\end{equation}
by integration on $[0,t]$. The entropy-energy inequality~\eqref{eq:En-En} follows by further integrating~\eqref{eq:I2} on $[0, + \infty)$ and using
$\displaystyle \lim_{t \to \infty} \ent{\mu}{P_t f} = 0$. 

Before making this argument rigorous we give a formal argument to~\eqref{eq:speed} at $t=0$, alternative to~\cite{ov00}. For simplicity, assume that $\mu=dx$ is 
the Riemannian measure and $(P_t)_{t \geq 0}$ is the heat semigroup associated with the Laplace-Beltrami operator $L = \Delta$.  Let $f$ be a probability density with respect to $dx$. First
$$
\partial_s P_{s\delta} f + \nabla\cdot(w_sP_{s\delta}f) =0,
$$
where $w_s=-\delta\nabla\log P_{s\delta}f$. Then one can check that at the first order in $\delta$, the couple $(P_{s\delta}f,w_s)_{s\in[0,1]}$ is optimal between $P_\delta f\mu$ and $f\mu$ in the Benamou-Brenier formulation (see~\cite[Chap. 7]{villani-book2}). Hence
$$
\frac{W_2^2(P_\delta f\mu,f\mu)}{\delta^2}=\int_0^1\int |\nabla\log P_{s\delta}f|^2P_{s\delta}fd\mu ds + o(1) \rightarrow I(f), \quad \delta \to 0.
$$

\medskip

%%%%%%%%%%%%%%%%
\begin{ethm}
In a REM space as in Section~\ref{sec-mms}, the contraction inequality~\eqref{eq-contraction-square-general-2} implies the entropy-energy inequality~\eqref{eq:En-En}. 
\end{ethm}
\begin{Proof}
 Let $f$ be a probability density with $f\mu \in \mathcal P_2(\X)$ and $I (f) < \infty,$ as we can assume. 
Recall that $( \X , d, \mu )$ is a $\mathsf{RCD} (R, \infty)$ space
under our assumption~\eqref{eq-contraction-square-general-2}. 
Thus, by \cite[Thm.~9.3 (i) and Thm.~8.5 (i)]{AGS13}, 
\begin{equation} \label{eq:speed2}
- \frac{d}{du}\ent{\mu}{P_uf} 
= 
I(P_uf) 
= 
\limsup_{\delta \downarrow 0} \frac{ W_2^2 ( P_{u+\delta} f \mu , P_u f \mu ) }{\delta^2}
\end{equation}
for a.e.~$u \in (0, + \infty)$. In particular,~\eqref{eq:speed} holds almost everywhere and, proceeding as above,
 \begin{equation}\label{eq:I}
I(P_tf)\leq e^{-2R(t-s)}I(P_sf)-\frac{2}{m}\int_s^t e^{-2R(t-u)}I(P_uf)^2du
\end{equation}
for any $t > s>0$ where~\eqref{eq:speed2} is valid.

\smallskip

We now prove that~\eqref{eq:I} holds for all $t>s \geq0$. 
%We have seen above that this formally implies~\eqref{eq:I2} by integration. In turn this concludes the argument.
For this, set $\psi (t) := e^{2Rt} I ( P_t f )$. 
Then $\psi$ is non-increasing on $[ 0, \infty )$ by a standard argument:
Indeed, by $CD ( R, \infty )$ with the self-improvement argument in \cite{Savare:2014jm}, we have 
$\sqrt{ \Gamma ( P_t f )} \le e^{-Rt} P_t ( \sqrt{ \Gamma (f) } )$ for all $t \geq 0$. It yields 
\[
\frac{ \Gamma ( P_t f ) }{ P_t f } 
\le e^{-2R(t-s)} \frac{ \left( P_{t-s} ( \sqrt{\Gamma (P_s f)} ) \right)^2 }{ P_{t-s} ( P_s f ) }
\le 
e^{-2R(t-s)}
P_{t-s} \left( 
    \frac{ \Gamma ( P_s f ) }{ P_s f }
\right).
\]
Thus the claim follows by integrating this inequality by $\mu$. 
Moreover $t \mapsto I ( P_t f )$ is lower semi-continuous (see e.g.~\cite[Lem.~4.10]{AGS13}). 
Thus $\psi$ is lower semi-continuous and non-increasing on $[0,\infty)$, so also right-continuous .
This implies that~\eqref{eq:I} holds for $t >s \geq 0$.

\smallskip
Let now $\delta > 0$. By dividing~\eqref{eq:I} by 
$e^{-2Rt}( \psi(t) + \delta ) ( \psi(s) + \delta )$,  for $t > s > 0$,
\begin{equation} \label{eq:pre_En-En1}
\frac{2}{m  ( \psi (s) + \delta ) ( \psi (t) + \delta ) } 
\int_s^t e^{-2Ru} \psi (u)^2 \, du 
\leq 
\frac{1}{ \psi (t) + \delta }
- 
\frac{1}{ \psi (s) + \delta } \cdot 
\end{equation}

 We claim
\begin{equation} \label{eq:pre_En-En2}
\frac{ 2 ( 1 - \delta ) }{m} 
\int_{0}^t e^{-2Ru} 
\left( 
    \frac{ \psi (u) }{ \psi (u) + \delta } 
\right)^2 
\, du 
\le 
\frac{1}{ \psi (t) + \delta }  
- 
\frac{1}{ \psi (0) + \delta }
\end{equation}
for any  $t \in [0, \infty)$. 
For the proof of the claim, we let $J$ be the subset of  $t \in [ 0 , \infty )$ 
satisfying~\eqref{eq:pre_En-En2}  and prove $J = [ 0 , \infty )$.
First, $0 \in J$ obviously holds and hence $J \neq \emptyset$. 
Second, if $t \in J$ and $t' \in ( t , \infty )$ with $t'- t$ sufficiently small, 
then $t' \in J$. Indeed, by the right continuity of $\psi$, 
we have 
$\psi (u) + \delta \geq ( 1 - \delta ) ( \psi (t) + \delta )$
for any $u> t$ being sufficiently close to $t$. 
 We take $t' > t$ so that 
this holds for all $u \in ( t , t' )$.
Thus~\eqref{eq:pre_En-En2} for this $t$,~\eqref{eq:pre_En-En1} and $\psi$ being non-increasing
yield 
\begin{align*}
\frac{ 2 ( 1 - \delta ) }{m}
\int_{0}^{t'} e^{-2Ru} 
&\left( 
    \frac{ \psi (u) }{ \psi (u) + \delta } 
\right)^2 
\, du 
\\
& \leq 
\frac{1}{ \psi (t) + \delta }  
- 
\frac{1}{ \psi (0) + \delta }
+ 
\frac{2}{m ( \psi (t) + \delta ) ( \psi (t') + \delta )} 
\int_{t}^{t'} e^{-2Ru} \psi (u)^2 \, du 
\\
& \le 
\frac{1}{ \psi (t') + \delta }  
- 
\frac{1}{ \psi (0) + \delta }
\end{align*}
and hence $t' \in J$. 
Third,  $J$ is closed under increasing sequences. 
That is, for  any bounded increasing sequence $( t_n )_{n \in \N}$ in $J$, then
$\displaystyle \lim_{n \to \infty} t_n  \in J$. 
This property follows from the fact that $\psi$ is  lower semi-continuous. 
Now these three properties imply  $J = [ 0 , \infty )$ and hence the claim holds. 

\smallskip

Finally we obtain~\eqref{eq:I2} for all $t \geq 0$
by taking $\delta \downarrow 0$ and rearranging terms 
in~\eqref{eq:pre_En-En2}.
 But 
\begin{equation}\label{HI}
\ent{\mu}{f} - \ent{\mu}{P_t f} = \int_0^t I ( P_s f ) ds 
\end{equation}
for all $t$ by \cite[Thms.~9.3 (i) and~8.5 (i)]{AGS13} again.
Hence integrating~\eqref{eq:I2}  in $t$ concludes the proof.~\end{Proof}

\bigskip

%%%
{\bf{A dimensional HWI type inequality}}

\medskip

For $R$ being $0$ or negative, no logarithmic Sobolev inequality for $\mu$ holds in general, and following~\cite{ov00} it can be replaced by a HWI interpolation inequality with an additional $W_2$ term : this is inequality giving an upper bound on the entropy $H$ in terms of the distance $W_2$ and the Fisher information $I$. As above, let us see how to derive a dimensional form of this inequality from the contraction property~\eqref{eq-contraction-sh}  in Theorem~\ref{thm-legros}. 

In a $REM$ space as in Section~\ref{sec-mms}, with a reference measure $\mu$ in $\mathcal P_2(\X)$, assume the contraction property~\eqref{eq-contraction-sh} with $R=0$. Let $f, g$ such that $f \mu, g \mu \in \mathcal P_2(\X)$, $I(f) < \infty$ and $g \mu$ has bounded support. Recall first that $( \X , d, \mu )$ is a $\mathsf{RCD} (0, \infty)$ space
under our assumption~\eqref{eq-contraction-sh}. In particular $I(P_t f) \leq I(f)$ for all $t \geq 0.$ Then \cite[Thm.~6.3]{AGMR} and the Cauchy-Schwarz inequality yield 

$$
\frac{1}{2} \frac{d}{dt}
W_2^2(P_t f \mu, g \mu )
\geq 
- W_2 ( P_t f \mu, g \mu ) \sqrt{I( P_tf)}
$$
for almost every $t>0.$ In particular
$$
\frac{1}{2}
W_2^2(P_t f \mu, g \mu )
-
\frac{1}{2}
W_2^2(f \mu, g \mu )
\geq - \int_0^t W_2 ( P_s f \mu, g \mu ) \sqrt{I( P_s f )} \, ds
\geq - \int_0^t W_2 ( P_s f \mu, g \mu ) \sqrt{I( f )} \, ds
$$
for all $t \geq 0$. 

If now $g$ converges to $1$ in such a way that $g \mu $ converges to $\mu$ in the $W_2$ distance, then using the triangular inequality
$$
\big\vert W_2(P_s f \mu, g \mu) - W_2(P_s f \mu, \mu) \big\vert \leq W_2(g \mu, \mu)
$$
for any $0 \leq s \leq t$ one can pass to the limit above, leading to
$$
\frac{1}{2}
W_2^2(P_t f \mu, \mu )
-
\frac{1}{2}
W_2^2(f \mu, \mu )
\geq - \int_0^t W_2 ( P_s f \mu, \mu ) \sqrt{I( f )} \, ds.
$$
 Now by~\eqref{eq-contraction-sh} the left-hand side is bounded from above by
 $$
 -4m \int_0^t \sinh^2 \left(\frac{ \ent{\mu}{P_s f} }{2m} \right) ds.
 $$
Finally $s \mapsto W_2(P_sf \mu, \mu)$ and $s \mapsto \ent{\mu}{P_s f}$ are continuous on $[0,t]$, so one can let $t$ go to $0$ and obtain
 %\[
%\frac{1}{2} \frac{d?{dt}\left( 
%    W_2^2(P_t f \mu, g \mu ) - W_2^2( f \mu, g \mu )
%\right)
%\geq 
%- W_2 ( f \mu, g \mu ) \sqrt{I(f)}.
%\]
%Hence, dividing~\eqref{eq-contraction-sh} by $t$  and letting $t$ go to 0  leads to 
\begin{equation} \label{eq-HWI1}
\sinh^2\left(\frac{ \ent{\mu}{f}}{2m} \right)\leq \frac{1}{4m}W_2(f\mu,\mu)\sqrt{I(f)}.
\end{equation}

A corresponding bound can also be derived for any $R$, in which the $s_{R/m}$ function appears.

\medskip

Here is a possible application of~\eqref{eq-HWI1}: in the above notation and assumptions (with $R=0$), there exists a positive numerical constant $C$ such that 
$$
\ent{\mu}{P_t f}  \leq \frac{m}{2} \max \Big\{ C ,  \log \frac{W_2^2(f \mu, \mu)}{mt} \Big\}, \quad t>0
$$
for all $f$ with $f \mu \in \mathcal P_2$. 
This bound is a consequence of~\eqref{HI}, \eqref{eq-HWI1} with $P_t f$ instead of $f$, the bounds $W_2(P_t f \mu, \mu) \leq W_2(f\mu, \mu)$ and $\sinh^4 (x) \geq e^{4x}/32$ for $x$ large enough.

%~\eqref{eq-cas-simple} with $s=0$ and again $P_t f$ instead of $f$.
For short time, this gives a regularization bound of the entropy as $m/2 \log (1/t)$, which is exactly the behaviour observed above for $R>0$, and also for the heat kernel on $\mathbb R^m$; it also improves on the corresponding bound $m \log (1/t)$ in~\cite[Prop. 2.17, (ii)]{EKS13}.

%Indeed, by the energy-dissipation identity \cite[(2.28)]{AGS13} for $\entf{\mu}$, 
%$\ent{\mu}{P_t f}$ and $\ent{\mu}{ P_t g }$ are continuous at $t = 0$, a.e. 

%For any $R$ we obtain
%$$
%2m\sinh^2\left(\frac{ \ent{\mu}{f} - \ent{\mu}{g} }{2m} \right)\leq \frac1{4} s_{\frac{R}{m}}(W_2(f\mu, g\mu))\sqrt{I(f)}-2Rs_{\frac{R}{m}}\left(\frac12 W_2(f\mu, g\mu)\right)^2
%$$

\medskip
\noindent
{\bf Acknowledgments.} This research was supported  by the French ANR-12-BS01-0019 STAB project 
and JSPS Grant-in-Aid for Young Scientist (A) 26707004.

%\footnotesize{\bibliographystyle{plain}
%\bibliography{biblio}
%}

\end{document}